\documentclass{sn-jnl}

\usepackage{graphicx}%
\usepackage{amsmath,amssymb,amsfonts}%
\usepackage[title]{appendix}%
\usepackage{xcolor}%
\usepackage{textcomp}%
\usepackage{manyfoot}%
\usepackage{booktabs}%
\usepackage{algorithm}%
\usepackage{algorithmicx}%
\usepackage{algpseudocode}%
\usepackage{listings}%
\usepackage[lined,boxed,commentsnumbered,algo2e]{algorithm2e}
\usepackage{geometry}
\usepackage{float}
\usepackage{url}

\begin{document}

\title[Article Title]{A Bottom-Up Approach to Optimizing the Solar Organic Rankine Cycle for Transactive Energy Trading}


\author[1,2]{\fnm{Silvia Anna} \sur{Cordieri}}

\author[3]{\fnm{Chiara} \sur{Bordin}}

\author[4]{\fnm{Sambeet} \sur{Mishra}}

\affil*[1]{Università di Bologna, silviaanna.cordieri2@unibo.it}

\affil[2]{Ricerca sul Sistema Energetico RSE S.p.A., silviaanna.cordieri@rse-web.it}

\affil[3]{The Artic University of Norway, chiara.bordin@uit.no}

\affil[4]{University of South-Eastern Norway, Sambeet.Mishra@usn.no}


\abstract{Solar Organic Rankine Cycle (ORC)-based power generation plants leverage solar irradiation to produce thermal energy, offering a highly compatible renewable technology due to the alignment between solar irradiation temperatures and ORC operating requirements. Their superior performance compared to steam Rankine cycles in small-scale applications makes them particularly relevant within the smart grid and microgrid contexts. This study explores the role of ORC in peer-to-peer (P2P) energy trading within renewable-based community microgrids, where consumers become prosumers, simultaneously producing and consuming energy while engaging in virtual trading at the distribution system level. Focusing on a microgrid integrating solar ORC with a storage system to meet consumer demand, the paper highlights the importance of combining these technologies with storage to enhance predictability and competitiveness with conventional energy plants, despite management challenges. A methodology based on operations research techniques is developed to optimize system performance. Furthermore, the impact of various technological parameters of the solar ORC on the system's performance is examined. The study concludes by assessing the value of solar ORC within the transactive energy trading framework across different configurations and scenarios. Results demonstrate an average 16\% reduction in operational costs, showcasing the benefits of implementing a predictable and manageable system in P2P transactive energy trading.}

\keywords{Solar organic Rankine cycles, peer-to-peer energy trading, transactive energy trading, storage systems, microgrids, operations research, energy systems modeling;}

\maketitle

\section{Introduction}
The new millennium has started with several innovations driven by the fast evolution of technologies in the energy sector \cite{CHICCO2009535}.  The consequences of climate change have posed significant challenges for governments worldwide. Despite the implementation of new directives to limit  those consequences, the most recent energy crisis highlighted the strong dependency that energy systems still have on conventional energy sources. Disruptive structural developments are still necessary to
deliver on the European Union’s COP21 commitments \cite{clean}, COP23 commitments \cite{CLEAN2} and UN sustainable goals \cite{CLEAN3}.

In this sense, the scientific community is working on solutions that minimize energy consumption and reduce environmental impact. Electricity production through the use of solar energy collection has been proven to be a viable option for green energy production \cite{YU2021114494}. Because of its abundance and availability, new solutions are continuously studied, to fully exploit its potential. Some studies are based on the idea of exploiting the full knowledge and experience gained over the last century on conventional power generation technologies and cycles but in a greener, decarbonized framework. When it comes to solar energy, a valid solution seems to be represented by Organic Rankine Cycles (ORCs). These cycles apply the same principles of a traditional steam Rankine cycle, replacing water as a working fluid with an organic fluid. Moreover, such cycles give the possibility to select the best working fluid and plant size depending on the available heat source. 

The possibility to select the best working fluid depending on the available heat source and the plant size results in multiple advantages: (i)
more efficient turbomachinery, (ii) limited vacuum at the condenser, and (iii) higher performance compared to both
steam Rankine cycles and gas cycles especially for heat sources lower than 400°C and power output lower than 20
MW \cite{TARTIERE20172}. Therefore, it seems perfectly suitable in a framework where a conventional heat source is substituted by a renewable one, i.e. solar energy. In fact, solar-driven technologies such as parabolic trough collectors can effectively produce heat at temperatures between 50 °C and 400 °C \cite{KALOGIROU2004231}. The framework consisting of an ORC driven by a solar heat source is referred to in the literature as the Solar Organic Rankine Cycle (Solar-ORC).


Energy informatics \cite{bordin2021educating} has become a key area of focus within energy research, driven by the need for intelligent and efficient energy systems that align with evolving sustainability and resilience goals. Smart grids, a significant component of this domain, integrate information and communication technologies with traditional power systems, enhancing their responsiveness, reliability, and efficiency. Recent advancements have been particularly focused on the modeling and simulation of smart energy and power systems, enabling predictive insights and dynamic management capabilities. These models aim to optimize energy flows, balance supply and demand, and facilitate the integration of renewable energy sources, thus supporting the transition to a more sustainable energy landscape \cite{bordin2020smart}
In such a system, vast numbers of devices, passively connected to the grid, will become actively involved in system-wide and local coordination tasks \cite{7452756}. In this context transactive energy trading emerges as a valid contender, to optimally coordinate such a complex scheme. The focus of this concept is mainly on the distribution level and its actors. Here, smart homes, buildings, and industrial sites engage in automated market trade with others at the distribution system level and with a two-way negotiation based on prices and energy quantities \cite{7452756}. In the future trend consumers become prosumers who can both produce and consume energy, but most importantly supply other consumers on a local level. This transactive energy trading among prosumers is called Peer-to-Peer (P2P) energy trading \cite{ZHANG2016147}. P2P is a decentralized form of transactive energy trading where prosumers are given the opportunity to engage without the need for an intermediary. This way, renewable energy integration is promoted either by investments in locally distributed energy resources made by prosumers or by encouraging consumers to purchase green energy locally, if they are incapable of investing in renewable energy sources. Although at the early stage, the P2P electricity trading without the need for utilities is expected to increase as the
awareness of the shared economy has grown and the microgrid has spread \cite{Park2017ComparativeRA}.  The main advantages of this system are: the power generation can be made meeting the requirements of the end users and the utilization of the resources can be optimized through the cooperative network between producers and consumers \cite{park2009technology}.\\
 A list of acronyms is provided in table \ref{nomenclature}.\\
 
\begin{table}
    \centering
    \begin{tabular}{l|l}
       \hline
       ORC  & Organic Rankine Cycle \\
        S-ORC & Solar-ORC \\
        TET & Transactive Energy Trading \\
        P2P & Peer-to-Peer \\
        MILP & Mixed Integer Linear Programming\\
        FPC & Flat-plate collector\\
        ETC&Evacuated tube collector\\
        CPC&Compound parabolic collector\\
        PTC&Parabolic through collector\\
        LFR&Linear Fresnel reflector\\
        \hline
    \end{tabular}
    \caption{List of acronyms}
    \label{nomenclature}
\end{table}

The objective of this paper is to investigate:
\begin{itemize}
    \item The compatibility between ORC and solar technology in very different locations weather-wise, Tromsø and Bologna.
    \item The potential that the Solar-ORC coupled with a storage system could have on a P2P transactive energy trading context.  Given the applications of this technology for reduced plant sizes, see Tartiere et al. \cite{TARTIERE20172}, it seems suited for the self-consumption requirements of a prosumer in such a trading context.
\end{itemize}
Moreover, we want to develop a tool that can optimize the management of the system we are considering. We do so, by means of operations research-based techniques. First, a MILP model for the operations scheduling of the Solar-ORC has been developed,  called the S-ORC model. Then, a MILP model has been implemented for the P2P Transactive Energy Trading between multiple prosumers in a local energy market where some Solar-ORCs are present as power generation plants owned by some prosumers, called the TET model.\\
 
The remainder of this paper is organized as follows. Section \ref{sec:Literature review} presents a literature review of the topics discussed in this paper, while Section \ref{sec:Novelty and key contribution} shows the novelty and key contributions of this work. Section \ref{sec:Technical notes, assumptions, and definitions} explains the main technical notes, assumption, and definitions, specifically in Section \ref{SM} there is a focus on the Solar-ORC, in Section \ref{TET} there is a focus on the transactive energy trading part, while in Section \ref{sec:solar collectors} there is a focus on the solar collector. Section \ref{sec:problem description} explains the problem we aim to solve. Section \ref{sec:Methodology} formally presents the S-ORC model and the TET model, discussed respectively in Section \ref{SM model} and \ref{TET model}, while Section \ref{sec:Model implementation} gives an insight on the model implementation. Section \ref{sec:Computational experiments} shows the computational experiments, specifically Section \ref{sec:Sensitivity analysis} presents a sensitivity analysis on the S-ORC model, while Section \ref{sec:Computational SM model} and Section \ref {sec:Computational TET model} discuss the computational experiments respectively on the S-ORC model and on the TET model. Section \ref{Discussions} contains further discussions and reflections on the computational results with a special focus on solar power plants in section \ref{sec:value solar}, while Section \ref{Future work} outlines possible research perspectives. Finally, Section  \ref{sec:Conclusions} draws conclusions.

\section{Literature review}\label{sec:Literature review}
This section, discusses the main contributions related
to transactive energy trading with a solar organic Rankine cycle problem. This analysis serves to contextualize the results that will subsequently be presented.\\

This paper combines the study of several topics, which in the past have been usually analyzed separately. Therefore, it seemed more functional to group all the contributions depending on the main topic they focus on, as one can observe in Table \ref{tab:State of Art}. The literature is classified based on the main modeling features of the treated problems. \\
\begin{table} 
    \centering 
    \resizebox{\textwidth}{!}{
    \begin{tabular}{|p{1,5cm}|p{1,7cm}|p{1,5cm}|p{1,5cm}|p{1,5cm}|p{1,5cm}|p{1,5cm}|p{1,7cm}|
}    \hline
    &Transactive energy trading&Peer-to-peer&Organic Rankine cycle&Solar-ORC&Energy storage&Simulation models&Prescriptive analytics\\
    \hline
      \cite{TARTIERE20172}&&&\checkmark&\checkmark&&&\\
     \hline
     \cite{TCHANCHE20092468}&&&\checkmark&\checkmark&&\checkmark&\\
     \hline
     \cite{YAMAGUCHI20062345}&&&\checkmark&\checkmark&&\checkmark&\\
     \hline
     \cite{7452756}&\checkmark&&&&&&\\
     \hline
     \cite{8966274}&\checkmark&\checkmark&&&&&\\
     \hline
     \cite{Park2017ComparativeRA}&\checkmark&\checkmark&&&&&\\
     \hline
     \cite{ZHANG20172563}&\checkmark&\checkmark&&&&&\\
     \hline  

     \cite{ZHAO2020100011}&&&\checkmark&&&\checkmark&\checkmark\\
     \hline

     \cite{PIEROBON2013538}&&&\checkmark&&&\checkmark&\checkmark\\
     \hline

    \cite{YU2021114494}&&&\checkmark&\checkmark&\checkmark&\checkmark&\\
     \hline
     \cite{ANEKE2016350}&&&&&\checkmark&&\\
     \hline

     \cite{CASATI2013205}&&&\checkmark&\checkmark&\checkmark&\checkmark&\\
     \hline

     \cite{BORDIN2019425}&&&&&\checkmark&\checkmark&\checkmark\\
     \hline

     \cite{BORDIN2019100824}&&&&&\checkmark&\checkmark&\checkmark\\
     \hline

     \cite{MISHRA2019696}&\checkmark&\checkmark&&&&\checkmark&\checkmark\\
     \hline

     \cite{MEHRPOOYA2017899}&&&\checkmark&\checkmark&\checkmark&\checkmark&\\
     \hline

    \cite{CHEN2022103824}&&&\checkmark&\checkmark&\checkmark&\checkmark&\\
     \hline

     \cite{WANG2020114327}&\checkmark&\checkmark&&&&\checkmark&\checkmark\\
     \hline

     \cite{ESMAT2021116123}&\checkmark&\checkmark&&&&\checkmark&\checkmark\\
     \hline 

     \cite{8282470}&\checkmark&\checkmark&&&&\checkmark&\checkmark\\
     \hline

     \cite{MAREFATI2021103488}&&&\checkmark&&\checkmark&\checkmark&\\
     \hline

     \cite{WANG2015401}&&&&&\checkmark&\checkmark&\checkmark\\
     \hline

\cite{MANFRIDA2016378}&&&\checkmark&&\checkmark&\checkmark&\\
     \hline

\cite{Fernandez}&\checkmark&\checkmark&&&&&\checkmark\\
     \hline 

\cite{Zhang}&\checkmark&\checkmark&&&&\checkmark&\\
     \hline     

     Our paper&\checkmark&\checkmark&\checkmark&\checkmark&\checkmark&\checkmark&\checkmark\\
     \hline
     
     \end{tabular}
     }
     \caption{Classification based on objective from the literature}
    \label{tab:State of Art}
\end{table} 
The discussion surrounding energy systems decentralization has drawn much attention among researchers to look into transactive energy trading, especially in a P2P framework. Kok et al. \cite{7452756} give an insight on the main coordination mechanisms of the smart grid vision, and on the role of transactive energy trading in this context. Zia et al. \cite{8966274} highlight potential reasons for avoiding the use of centralized microgrid transactive energy systems, and discuss existing architectures for a decentralized transactive energy system. Park et al. \cite{Park2017ComparativeRA} 
provide a comprehensive review of the design of peer-to-peer markets, as well as their challenges and opportunities, while Zhang et al. \cite{ZHANG20172563} discuss existing P2P projects. The scientific community's significant interest in P2P energy trading has produced different strategies to tackle this problem. Esmat et al. \cite{ESMAT2021116123} propose a platform for decentralized P2P trading based on two key layers. A market layer features a short-term multi-staged multi-period market with a uniform pricing mechanism. Then a blockchain layer offers a high level of automation, security, and fast real-time settlements through smart contract implementation. Mishra et al. \cite{MISHRA2019696} develop a multi-agent approach, where a math-heuristic model is used in the context of a decentralized power distribution system. In their work Wang et al. \cite{WANG2020114327} present a method based on the double auction market. Here each prosumer firstly dispatches its flexible energy resources with the objective of minimum cost, then the coordination of energy resources among diversified prosumers can be achieved with the aid of P2P energy transactions. Finally, Khorasany et al. \cite{8282470} implement a platform, where prosumers with excess energy and consumers communicate with each other to maximize their welfare. A double auction with an average mechanism is applied to determine the allocation and price of energy. Fernandez et al. \cite{Fernandez} developed a Community Energy Management System (CEMS)
is presented in this paper to facilitate local P2P trading among
consumers based on bi-level optimization
to maximize the utility of all involved parties. Zhang et al. \cite{Zhang} implemented a model for P2P energy trading and an associated bidding system for the P2P energy trading among consumers
and prosumers in a grid-connected Microgrid, consisting of a four-layer
system architecture. \\

The role of solar-ORCs has been widely discussed in the literature.  Zhao et al. \cite{ZHAO2020100011} provide a detailed literature review on each design procedure of ORCs using artificial intelligence algorithms. A comprehensive view of the ORC market is given by Astolfi et al. \cite{TARTIERE20172}, explaining the main ongoing applications and the role of solar-ORCs. Pierobon et al. \cite{PIEROBON2013538} show a multi-objective optimization with a genetic algorithm for the optimal design of ORCs. Other papers study specifically on solar-ORCs, focusing on different aspects. Tchanche et al. Some works such as \cite{TCHANCHE20092468}- \cite{YAMAGUCHI20062345} investigate the impact of different organic working fluids on the plant's overall performance, while Chen et al. \cite{CHEN2022103824} introduce and evaluate using Aspen-HYSYS and MATLAB software, a Solar-ORC configuration where solar energy plays a key role in the production of energy and hydrogen fuel. Here the ORC is fed by a solar farm based on the parabolic trough solar collector (PTSC), and then a fraction of the electrical energy obtained is fed into an alkaline electrolyzer (AEL) to produce hydrogen fuel. Mehrpooy et al. \cite{MEHRPOOYA2017899} concentrate on the design optimization of the Solar-ORC, which is evaluated through a thermoeconomic performance. The optimal point was selected using the TOPSIS decision-making technique among the Pareto frontier of the genetic algorithm. Finally, Yu et al. \cite{YU2021114494} implement a simulation-based optimization model in Aspen HYSYS to optimize both the design and operation of a Solar-ORC.\\

One major challenge facing a solar-driven energy source such as Solar-ORC, is the intermittency which makes it unreliable for a steady energy supply. Through the energy storage concept, these renewable resources can be made to be reliable and steady energy sources \cite{ANEKE2016350}. The coupling of energy generation and storage has become a trend nowadays in the scientific community. Casati et al. \cite{CASATI2013205} study 
the role of thermal energy storage for a Solar-ORC. Manfrida et al. \cite{MANFRIDA2016378} focus on a robust mathematical model of a Latent Heat Storage (LHS) system constituted by a storage tank containing Phase Change Material spheres. The model is simulated under dynamic (time-varying) solar radiation conditions with the software TRNSYS.
Marefati et al. \cite{MAREFATI2021103488} present the performance study of a Pumped-Hydro and Compressed-Air storage system, coupled with an organic Rankine cycle (ORC).
Wang et al. \cite{WANG2015401} implement an LP optimization model for a combined heat and power (CHP) based DH system with RES and energy storage system (ESS). Finally, some papers such as \cite{BORDIN2019100824}-\cite{BORDIN2019425} and \cite{bordin2017linear} (the latter further expanded in the work \cite{bordin2015mathematical}) propose optimization models that include battery degradation, and highlight its impact on having realistic performance of such systems.\\

While many works in the literature address some of the topics covered in this paper, the majority do it separately, featuring just some of them.\\  
Therefore there exists a research gap in the form of: 

\begin{itemize}
 
\item Technological representation of solar ORC in ways suitable for inclusion within mathematical optimization models for operational planning of energy systems.  

\item A practical understanding of the value of solar ORC in peer-to-peer interaction at the microgrid level.
\end{itemize}

To the best of our knowledge, the problem we introduce in this paper is the first to simultaneously feature a MILP model for transactive energy trading in a P2P context for a solar-ORC coupled with a storage system.

\subsection{Novelty and key contribution}\label{sec:Novelty and key contribution}

The main contributions of this paper can be divided into two categories: a methodological contribution and an analytical contribution.\\

From a methodology point of view, we propose a MILP model for solar-ORCs coupled with a storage system that includes technological details. More specifically, the model considers detailed energy balances for the components of the cycle. It also contains thermodynamic properties of the fluid to see how different working fluids impact the performance of the plant. Moreover, battery degradation is also included, to optimize battery usage.
This
model can be inserted into a wider optimization model for P2P transactive energy trading. Both models can be used as stand-alone models or can be easily included in large open-source energy system model. Traditional energy systems models available in the literature, especially the largest ones, are usually technology agnostic and, thus, do not contain a detailed description of the technologies involved. In fact, technologies are usually treated as black boxes without considering technological features. This simplification may be functional for certain problems and may be more competitive from a computational time point of view. However, including technological details can be the key to more realistic implementations and results. Moreover, having such details can open possibilities for new users, that specifically request such information.\\

The analytical contribution given by this paper is represented by an extensive sensitivity analysis. First, different frameworks are tested for the Solar-ORC with several working fluids, plant sizes, and solar collector's technologies. We implement this analysis to understand the general value of the Solar-ORC. Then the Solar-ORC is introduced in a P2P transactive energy trading framework to perform more analyses. The aim is to evaluate the value of such a system in this context and to better understand the benefit that this plant system could have on the community. Different scenarios, i.e. different energy communities are considered. We create instances that represent domestic and industrial users who can also be prosumers, thus satisfying their own demand loads. Finally, all the tests are repeated for different cities and different seasons of the year, to understand how the performance and value of the system may be affected.\\

\section{Technical notes, assumptions, and definitions}\label{sec:Technical notes, assumptions, and definitions}

This section, gives an insight of which are the main technical aspects and assumptions concerning this work. Such insight is substantial for the reader to better understand the models that we propose later on. More specifically in section \ref{SM} the insight is referred to the Solar-Organic Rankine cycle, while section \ref{TET} is focused on the transactive energy trading part, finally section \ref{sec:solar collectors} discusses the solar collector specifics.

\subsection{Solar-Organic Rankine cycle}\label{SM}
\begin{figure} 
    \centering
    \includegraphics[width=\textwidth,keepaspectratio]{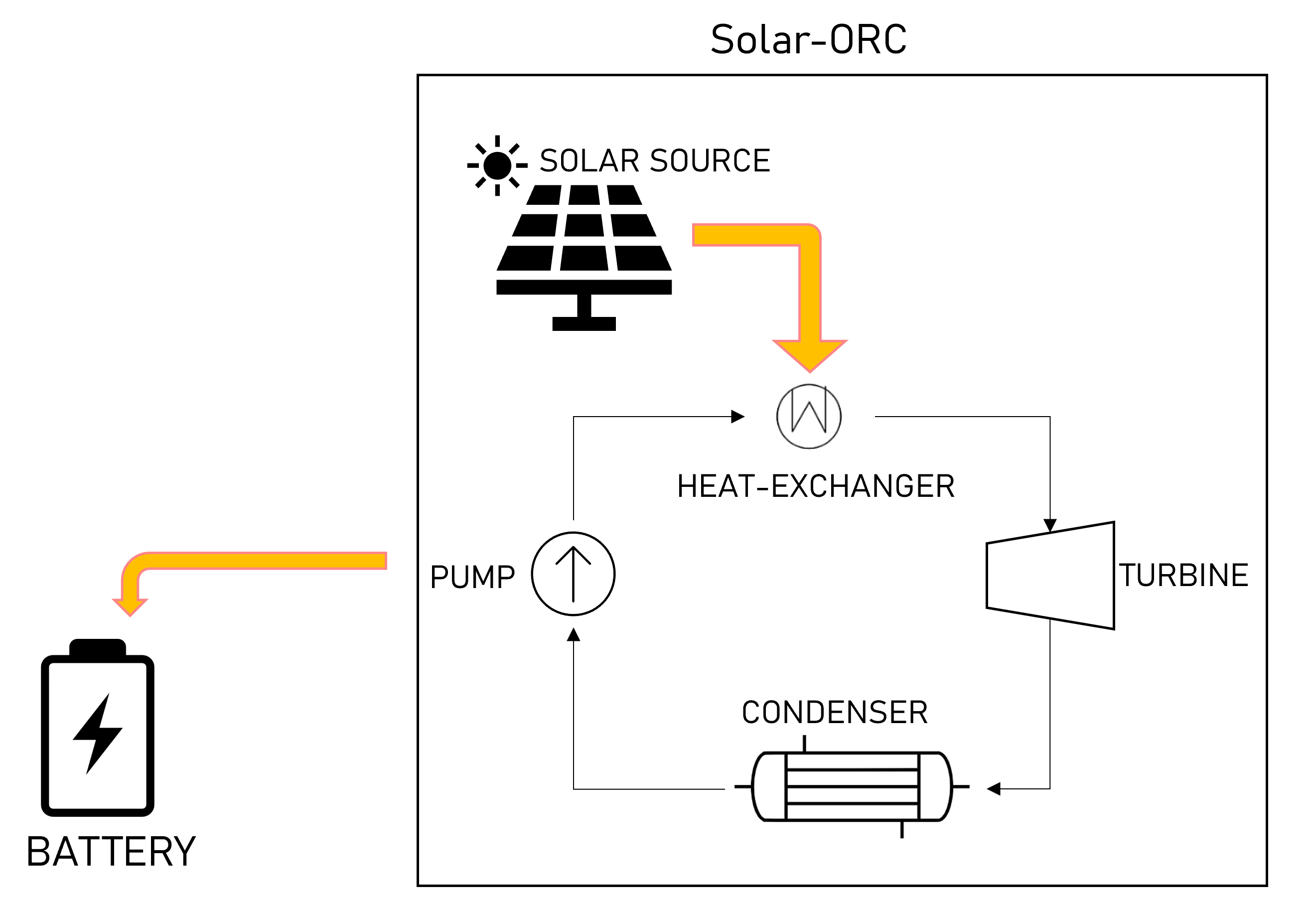}
    \caption{Scheme of the Solar-ORC}
    \label{fig:S-ORC}
\end{figure}

Solar-Organic Rankine cycles are characterized by using the sun as a source of thermal energy. In fact, a solar collector acts as an evaporator to heat the working fluid of the Rankine cycle. The use of solar irradiation for driving an ORC is a promising
renewable energy-based technology due to the high compatibility between the operating temperatures of solar
thermal collector technologies and the temperature needs of the cycle \cite{LONI2021111410}. In fact, organic Rankine cycles usually operate at
temperatures of up to 400 °C or 500 °C, which is perfectly compatible with thermal energy available from solar-based technologies.\\
The Solar-ORC scheme is depicted in Figure \ref{fig:S-ORC}. The ORC sub-system consists of a pump, an evaporator, a turbine and a condenser. The organic working fluid is pumped from condensation pressure to evaporation pressure. After pumping, the organic working fluid is vaporized and superheated in the evaporator, using thermal energy supplied by solar panels. Next, the high temperature and high-pressure vapor is expanded through the turbine to generate power. Finally, the working fluid is condensed in the condenser.\\
The Solar-ORC is connected to an electricity storage system consisting of a battery. The battery is used to store electricity, whenever there's an overproduction of the Solar-ORC, which can be withdrawn whenever it is needed.

The scope of our research work surrounds the economic operational optimization of the plant's scheduling. Thus, we include in our optimization process the components of the cycle, that are directly connected to the net power output of the plant, i.e. 
the pump, the turbine, and the heat exchanger connected to the solar panels.\\

The regulation of the steam turbines is used at a constant velocity to adapt the power of the turbine. We assume to apply lamination as a regulation policy of the cycle. In
this regulation the process is
at constant enthalpy. By closing a valve, thus reducing the section area, at the entrance of the turbine the pressure of the steam is reduced, while the entropy rises.
As the valve is closed, the constant enthalpy process occurs through the valve
with an increase in entropy and a decrease in the availability of energy per unit of mass
flow rate.\\ 

The thermodynamic properties of the working fluid are fundamental to determining the economics of an ORC. A bad choice could lead to a low efficient and expensive plant \cite{TCHANCHE20092468}. Thus, we include in our model a more detailed calculation of the mass flow rate of the working fluid, directly dependent on the type of organic fluid used. We choose to use density as the parameter representative of the thermodynamic properties of the organic working fluid.\\
9 different types of organic fluid, already well-known in the literature, are presented. These fluids are considered valid candidates. However, the methodology that we propose in Section \ref{sec:Methodology} can be applied to any type of organic fluid.\\

Energy storage is an essential link in the energy supply chain \cite{ANEKE2016350}. This is enhanced when it comes to most renewable energy resources, especially solar and wind. As a matter of fact, they occur intermittently, which makes them unreliable for a steady energy supply. If coupled with energy storage technologies, these renewable resources can increase their reliability.\\
The battery system technology is the most widespread energy storage device for power system application \cite{DIVYA2009511}. They seem to be a commonly applied solution lately to deal with renewable energy sources' instability. Nonetheless, the gain obtained in stabilizing the system may not be proportional to the increase in costs. In fact, the installation of batteries does not always automatically reduce the cost enough to pay for the installation \cite{BORDIN2019425}. Therefore, including an optimization of the storage system from a usage point of view, may be crucial to contain economic losses. In fact, the lifetime of a battery is highly influenced by the way it is operated, and by deterioration. Bad handling could result in more frequent substitutions of the battery, thus in higher costs.\\
The parameter that measures the life of a battery is called lifetime throughput. It 
defines the total dischargeable amount of energy in kWh before it is no longer able to deliver energy, enough to satisfy the load requirements of the system \cite{BORDIN2019100824}. The residual number of cycles to failure is inversely proportional to the depth of discharge. Deeper discharge results in a lower number of related cycles to failure.\\
Another important parameter is the state-of-health of a battery. This is a percentage of the battery capacity available when fully charged relative to its rated capacity. The state-of-health accounts for battery aging.\\
In our paper, we include both the lifetime throughput and the battery fade due to aging in our methodology. This way, the optimization process will avoid a non-economically optimal use of the storage system. \\

\subsection{Transactive energy trading}\label{TET}
Transactive energy trading emerges as a valid option among smart grid handling tools. The concept of having local actors that handle the grid on a distribution level opens the opportunity for consumers to be more active and involved. Considering that self-consumption has become greatly widespread, thanks also to incentives given in the last decades by governments, the role of consumers has substantially changed. The so-called prosumers are now important actors, that can no longer be considered as passive entities. Especially, when the stability of the grid is involved, it becomes even more clear that new management options are necessary to deal with these deep changes in the grid's framework.\\
From the prosumers' point of view, especially those less experienced, having an
automated market
trade can be crucial. This way they can concentrate on the optimization of their own profit. In a self-consumption framework, the prosumers invest to fulfill their own consumption, usually using a renewable energy source. Thus, the profit comes mainly from handling his own demand. The possibility of selling overproduction to other prosumers, or buying when there's a lack of production, is a plus. In view of this concept, the goal of our work
in the trading phase is to only optimize the under/over supply of electricity among the micro-grids of
our system. As we will discuss more in detail later on, this also will give us a great advantage from a computational point of view.\\

In this paper, we concentrate
on short (hour)-medium (week) operational planning. We consider this to be more meaningful for the problem we are inspecting. Therefore, we concentrate on a weekly time range to lead computational experiments.\\

\begin{figure} 
    \centering
    \includegraphics[width=\textwidth,keepaspectratio]{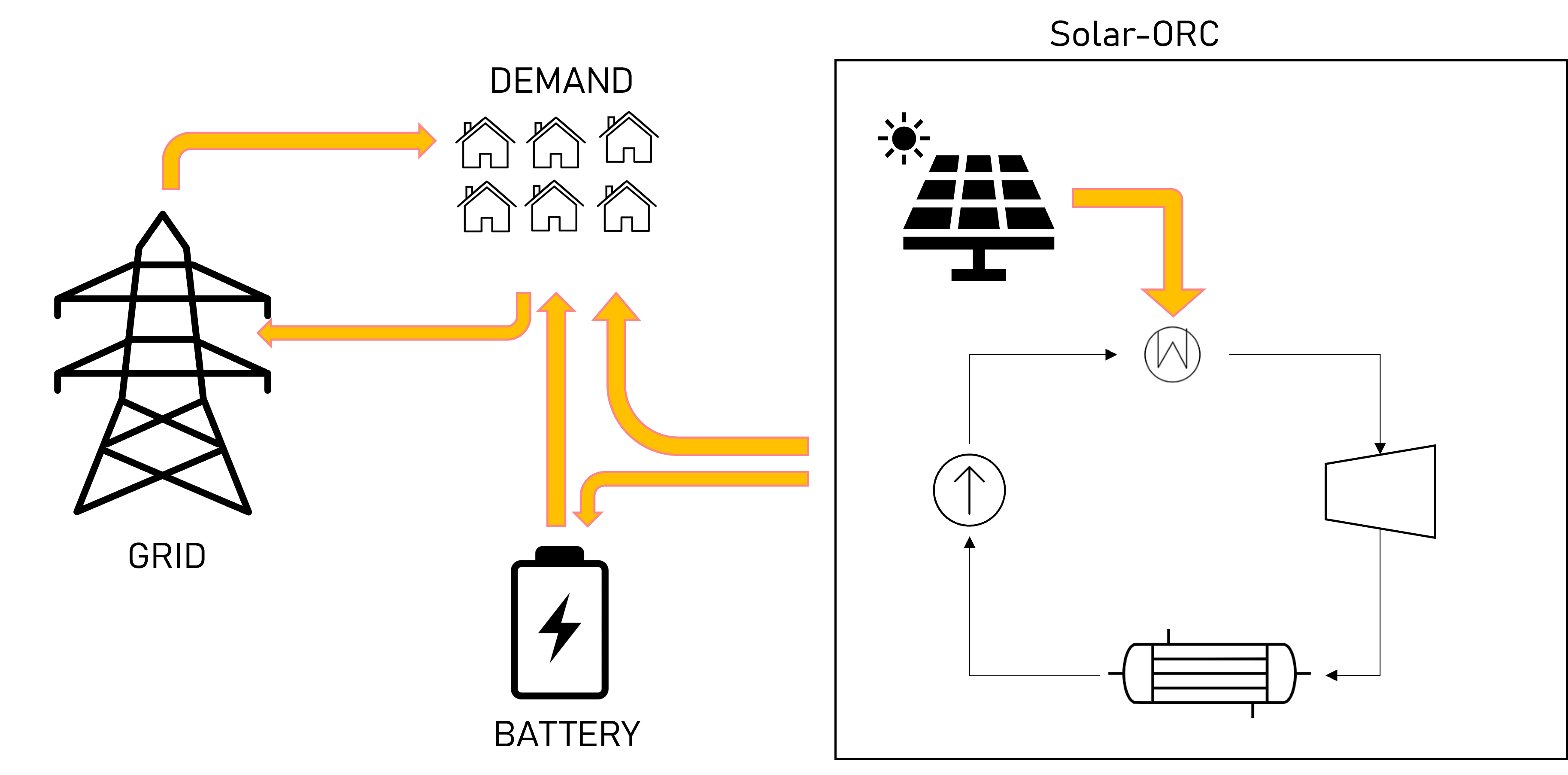}
    \caption{Scheme of the system considered}
    \label{fig:TET scheme}
\end{figure}

Figure \ref{fig:TET scheme} shows the scheme of the single microgrid that we want to optimize. The Solar-ORC is coupled with a battery, to fulfill the prosumer's demand. The grid is used to balance over/under supply by the Solar-ORC.

\begin{figure} 
    \centering
    \includegraphics[width=0.9\textwidth,keepaspectratio]{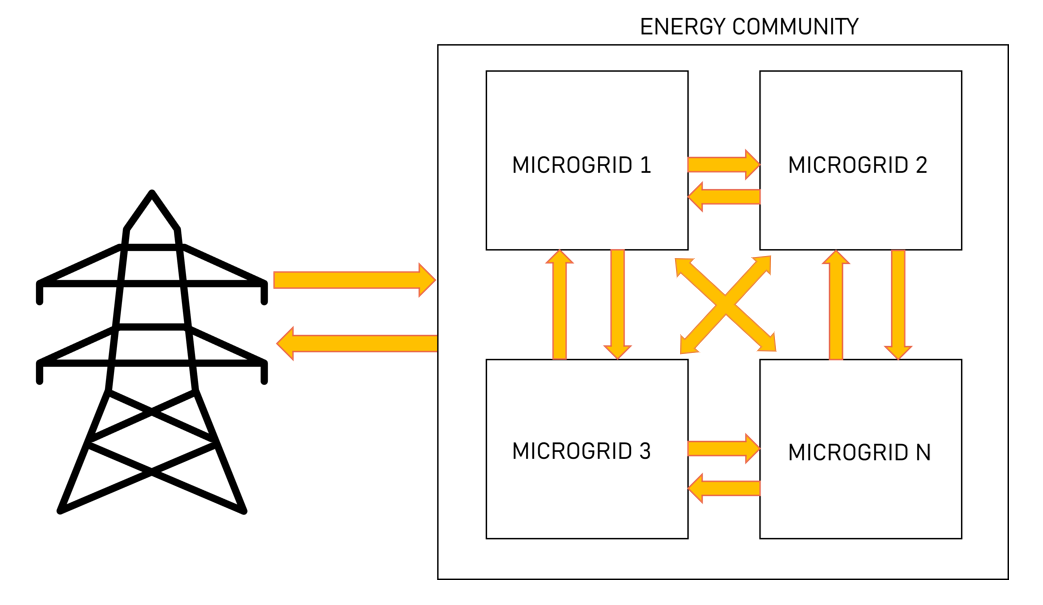}
    \caption{Scheme of the system used for the P2P trading}
    \label{fig:TET scheme 1}
\end{figure}

Figure \ref{fig:TET scheme 1} shows the total P2P transactive energy trading system that we want to optimize in our work. Here the microgrids that represent the single prosumers are able to interact with each other. The grid is used to balance the total over/under supply of the prosumers.

\subsection{Design and operational parameters
of the solar collector}\label{sec:solar collectors}
The primary energy source of the Solar-ORC is solar energy, which is provided to the cycle through a solar collector. The solar energy provided by the solar collector depends on the efficiency of the collector chosen, the area of the collector, and the beam irradiation. An ORC engine can be coupled with various solar collectors because the power cycle can operate in a large range of heat source temperatures \cite{LONI2021111410}. Each collector technology has a different efficiency, that determines the actual primary energy available to the cycle. According to Loni et al. \cite{LONI2021111410}, there are three solar collectors categories that are particularly suitable for Solar-ORC applications: flat-plate collectors, solar concentrators, and tracking solar concentrators.\\

\begin{figure} 
    \centering
    \includegraphics{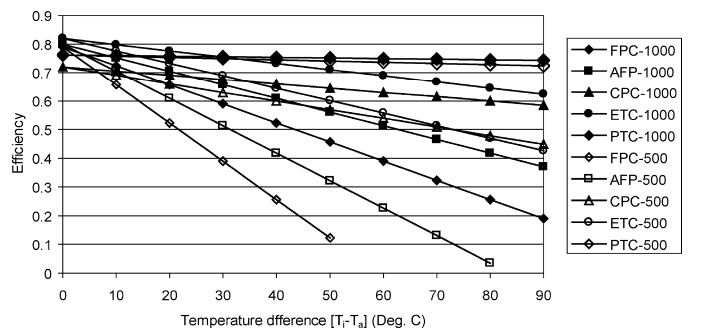}
    \caption{Comparison of the efficiency of various collectors presented by Kalogirou et al. \cite{KALOGIROU2004231}.}
    \label{fig:efficiencies}
\end{figure}

In this paper five different technologies are tested, three flat-plate collectors (Flat-plate collector (FPC), Evacuated tube collector (ETC) and Compound parabolic collector (CPC)) and two solar concentrators (Parabolic through collector (PTC), Linear Fresnel reflector (LFR)). The objective is to analyze different solar collector and ORC couplings, to understand the impact on the system considered. The main parameter that we will use to determine the impact is the solar collector's efficiency. Therefore, this parameter is considered as input data in our model as will be better explained in section \ref{sec:Methodology}.\\
ASHRAE Standard 93:1986 \cite{ASHRAE} is undoubtedly the one most often used for testing the thermal performance of collectors. It gives information on
testing solar energy collectors. The data can be used
to predict performance in any location and under any
weather conditions where load, weather, and insolation are
known. There are studies in the literature that based on ASHRAE Standard 93:1986, derive efficiency curves for different solar collectors' technologies. An example is depicted in Figure \ref{fig:efficiencies}. The solar collectors' efficiency values used are taken from the literature \cite{KALOGIROU2004231}.

\section{Problem description}\label{sec:problem description}
This section shows the problems we aim to solve in this paper. First, we want to optimize the operational scheduling of the system depicted in Figure \ref{fig:TET scheme}. The idea is to schedule the functioning of a single Solar-ORC, coupled with a battery as a storage system to fulfill the final consumer's demand. The battery is subject to degradation, with respect to the charging/discharging cycles. More specifically, a non-efficient usage diminishes its maximum capacity. The grid is used to balance over/under production of the Solar-ORC.\\
Then, the system depicted in Figure \ref{fig:TET scheme 1} is considered. Here the first scheme is included in a local community. The local community consists of several prosumers trading energy with each other in a P2P way. The main objectives of this second problem are overall operational cost minimization while prioritizing self-consumption for each single prosumer.\\ 
In both problems, we consider a time horizon of one week, divided into hourly time intervals.\\

\section{Methodology}\label{sec:Methodology}

This section presents the methodology used and the two MILP models that is used to optimize the problem previously presented. Section \ref{sec:Model implementation} explains the implementation process considered. Section \ref{SM model} presents an optimization model for the management of S-ORCs coupled with a storage system, called the  S-ORC model. Section \ref{TET model} we presents an optimization model for Transactive Energy Trading, called the TET model, in a P2P context for a Solar-ORC coupled with a storage system.
Table \ref{tab:sets} summarises all the sets parameters and variables used in the S-ORC model and in the TET model.\\
\begin{table}  
    \centering
    \begin{tabular}{l l}
    \hline\\ 
    Sets&\\
    \hline\\
    $t\in  \{0,..,T\}$ & Set of time intervals \\
    $i \in N$ & Set of microgrids\\
    \hline\\
    Parameters&\\
    \hline\\
    $\eta_{th}$ & Efficiency of the heat exchanger\\
    $\eta_{b}$ & Charging/discharging efficiency of the storage system\\
    $\eta_{I}$ & Efficiency of the Solar-ORC\\
    $\eta_{solar}$ & Efficiency of the panels\\    
    $c_p$ & Cost of production\\
    $c_{t,ij}^t$ & Transmission cost from $i \in N$ to $j \in N$\\ 
    $c_b$ & Cost of charging/discharging\\
    $x_{min}$ & Minimum power boundary of the ORC\\
    $x_{max}$ & Maximum power boundary of the ORC\\   
    $z_{min}$ & Minimum power boundary of the pump\\
    $z_{max}$ & Maximum power boundary of the pump\\    
    $g_{min}$ & Minimum boundary of electricity that can be injected in the grid\\
    $g_{max}$ & Maximum boundary of electricity that can be injected in the grid\\
    $b_{min}$ & Minimum power boundary of the storage system\\
    $b_{max}$ & Maximum power boundary of the storage system\\    
    $D^t$ & Demand of each $t \in \{0,..,T\}$\\
    $v$ & Velocity of the working fluid\\
    $A_{solar}$ & Area of the  solar panels\\
    $I^t_{solar}$ & Beam irradiation for every $t \in T$\\
    $f_{ij}^{max}$ & Maximum flux  of electricity limit from prosumer $i\in I$ to prosumer $j \in J$\\
    $f_{ij}^{min}$ & Minimum flux of electricity limit from prosumer $i\in I$ to prosumer $j \in J$\\ 
    $\rho$ & Density of the working fluid\\    
    $\Delta h_P$ & Enthalpy difference in the pump\\
    $\Delta h_T$ & Enthalpy difference in the turbine\\
    $B_{fade}$ & Battery fade in efficiency due to aging\\
    $B_{throughput}$ & Battery throughput\\
     \hline\\
    Variables&\\
    \hline\\
    $g^t$ & Injection of electricity in $t\in \{0,..,T\}$\\
    $q_{in}^t$ & Thermal power coming from the heat exchanger for every $t\in \{0,..,T\}$ \\
    $z^t$ & Power consumed by auxiliaries and pump in the ORC\\
    $b^t$ & Battery level for every $t \in \{0,..,T\}$ \\
    $b^t_{in}$ & Power flow entering the battery for every $t \in \{0,..,T\}$ \\  
    $b^t_{out}$ & Power flow injected in the grid for every $t \in \{0,..,T\}$ \\
    $y^t_{in}$ & Binary variable connected to the power flow entering the battery for every $t \in \{0,..,T\}$ \\ 
    $y^t_{out}$ &  Binary variable connected to the power flow injected in the grid for every $t \in \{0,..,T\}$ \\ 
    $b^t_{max}$ & Maximum capacity of the storage system at time $t \in \{0,..,T\}$ \\
    $d^t$ & Storage system degradation\\
    $q_{solar}^t$ & Solar power injected in the heat exchanger for $t \in \{0,..,T\}$\\
    $x^t$ & Power produced by the 
    Organic Rankine cycle for $t\in \{0,..,T\}$\\
    $m_{ORC}^t$ & Mass flow rate of the Organic Rankine Cycle for $t\in \{0,..,T\}$\\
    $A^t$ & Section area traversed by the mass flow rate for $t\in \{0,..,T\}$ \\    
    $f_{ij}^t$ & Flux of electricity sold from microgrid $i\in I$ to consumer $j \in J$ every $t \in $\\
    $e_{in}^t$ & Electricity taken from the grid by the single microgrid every for $t\in \{0,..,T\}$\\
    $e_{out}^t$ & Electricity sold to the grid by the single microgrid every $t\in \{0,..,T\}$ \\   
    $h_{in}^t$ & Electricity taken from the grid to balance the whole system every $t\in \{0,..,T\}$\\
    $h_{out}^t$ & Electricity sold to the grid to balance the whole system every $t\in \{0,..,T\}$\\    
    
    \hline
    \end{tabular}
    \caption{Sets, parameters and variables used in the S-ORC model and in the TET model}
    \label{tab:sets}
\end{table}

\subsection{Model implementation}\label{sec:Model implementation}
To optimize the two problems, we implement algorithms based on two mathematical models, that are better explained later in section \ref{SM model} and section \ref{TET model}. The mathematical models are developed following Operations Research theory, specifically Mixed Integer Linear Programming (MILP) theory. The idea is that each value that needs to be optimized in the system, is represented by a variable. The variables are then used to formulate the constraints and objective function of the problem. 
Once formulated, the mathematical models can be easily implemented into a script using different programming languages (e.g. Julia \cite{JULIA}) and solved to optimality through a solver. The solver is a commercial software, that finds the optimal solution using a combination of mathematical programming techniques. There are several solvers that have already been developed and are currently available to final users. The specifics of the solver used in this work are presented in section \ref{sec:Computational experiments}.

\subsection{S-ORC model}\label{SM model} 
In the following, an optimization model for the management of S-ORCs coupled with a storage system is presented, called the S-ORC model.\\
The planning horizon is divided into $t \in T$ hourly intervals. Note that this can be adjusted to meet the needs of the application of the model. For each time interval, the final consumer's demand is given $D^t$. Several parameters are related to the thermodynamic characteristics of the system. The efficiency of the Solar-ORC is denoted by parameter $\eta_{I}$, while $\eta_{th}$ is the efficiency of the heat exchanger.  The maximum and minimum electric energy production limits of the Solar-ORC are $x_{max}$ and  $x_{min}$, whereas $g_{max}$ and  $g_{min}$ are used to limit the electricity available to be used. In this sense, the grid is used in this model to balance the system in case of over/under production. Thus, the overproduction can be injected into the grid and vice versa underproduction can be withdrawn from the grid. The velocity $v$ and density $\rho$ of the working fluid are given, as the specific heating value $C_{P,ORC}$. The Solar-ORC is characterized by a temperature difference inside the pump $\Delta T_{P}$ and inside the turbine $\Delta T_T$, whereas the efficiencies of the pump and the turbine are denoted as $\eta_P$ and $\eta_T$.\\
The solar part of the S-ORC plant is defined by the total area $A_{solar}$ and the efficiency of the solar panel $\eta_{solar}$. For every time interval $t \in T$ the beam irradiation $I_{solar}^t$ is given.\\
The storage system is defined by the charging/discharging efficiency $\eta_{b}$, and by a maximum and minimum limit the capacity, respectively $b_{max}$ and  $b_{min}$. We consider a battery as a storage system. To better evaluate the actual capacity of the battery, we introduce two parameters, the lifetime throughput of the battery $B_{throughput}$, and the battery fade $B_{fade}$. The lifetime throughput measures the life of a battery. It defines the total amount of energy in kWh that can be discharged before it cannot satisfy the load requirements of the system. Additionally, the battery fade is used to calculate the loss of capacity as the battery ages. The capacity of the battery will not drop more than a certain percentage, i.e. $B_{fade}$ as long as the total energy drawn is kept within the lifetime throughput \cite{BORDIN2019100824}.\\
Every kWh of electricity produced by the Solar-ORC has a cost $c_p$ while every kWh of electricity stored has a storage cost of $c_b$.\\

The objective of the model is to minimize the total costs of production, given by the cost of production and the cost of storage. The electricity produced by the turbine of the Solar-ORC every time interval $t \in T$ is measured by variable $x^t$, while variables $b^t_{in}$ and $b^t_{out}$ measure the electricity respectively charged or discharged every $t \in T$.

\begin{align}\label{objective function}
 [\text{S-ORC model}]\    \text{minimize}   \sum_{t \in T} c_{p}x^t+c_{b}(b_{in}^t+b_{out}^t)
\end{align}
Constraints (\ref{grid 1})-(\ref{grid 3}) define variable $g^t$, the electricity available to be used every $t \in T$. Specifically, constraints (\ref{grid 1}) limit its capacity, while constraints (\ref{grid 2})-(\ref{grid 3}) are energy balances on the system, in agreement with the literature \cite{Cordieri}. Furthermore, variables $e_{in}^t$ and $e_{out}^t$ indicate the amount of electricity withdrawn or injected in the grid. 

\begin{align}
   \label{grid 1}
   g_{min}\leq g^t\leq g_{max}\  \forall t \in T\\
   \label{grid 2}
   g^t = x^t-z^t-\eta_{b}b^t_{in}+\frac{b^t_{out}}{\eta_{b}}\   \forall t \in T  \\
     \label{grid 3}
        g^t + e_{in}^t \geq D^t + e_{out}^t \ \forall t \in T 
\end{align}
Constraints (\ref{grid in 4})-(\ref{pump power}) describe the energy balances to define the actual electricity production of the Solar-ORC. More precisely constraints (\ref{grid in 4}) connect the net energy produced, given by the subtraction of the energy produced by the turbine $x^t$ and, the one consumed by the pump, $z^t$, to the thermal energy coming from the heat exchanger for every $t \in T$. While constraints (\ref{power turbine}) and (\ref{pump power}) measure the energy produced by the turbine and consumed by the pump respectively. Both these energies are regulated by the mass flow rate of the working fluid represented by variable $m^t_{ORC}$. The value of energy produced by the turbine and consumed by the pump is limited by constraints (\ref{x limit})-(\ref{z limit}) respectively.
\begin{align}
    \label{grid in 4}
  x^t-z^t = \eta_{I}q_{in}^t\  \forall t \in T \\
    \label{power turbine}
  x^t=m^t_{ORC}\Delta h_{T} \ \forall t \in T\\
    \label{pump power}
  z^t=m^t_{ORC}\Delta h_{P} \ \forall t \in T\\
    \label{x limit}
  x_{min}\leq x^t\leq x_{max} \ \forall t \in T\\
    \label{z limit}
  z_{min}\leq z^t\leq z_{max} \ \forall t \in T
\end{align}  
The mass flow rate of the working fluid for every $t \in T$ is calculated through constraints (\ref{mass flow rate}). Here variable $A^t$ represents the actual section of the pipes, regulated every $t \in T$. The pipe section is regulated following the lamination concept \cite{Lozza}.
\begin{align}
      \label{mass flow rate}
  m^t_{ORC}= \rho A^t v ,\ \forall t \in T
\end{align}
The solar part of the S-ORC plant is managed through constraints (\ref{grid in 5}) and (\ref{solar heat}). In fact, constraints (\ref{grid in 5}) link the thermal energy provided by the solar panels to the thermal energy available at the heat exchanger of the Solar-ORC, whereas constraints (\ref{solar heat}) compute the thermal energy provided by the solar panels with respect to the beam radiation. 
\begin{align}  
  \label{grid in 5}
  q_{in}^t \leq
  \eta_{th}q_{solar}^t\  \forall t \in T \\
  \label{solar heat}  q_{solar}^t=\eta_{solar}A_{solar}I_{solar}^t \ \forall t \in T
 \end{align}
The battery management is provided by constraints (\ref{battery1})-(\ref{battery in}). More specifically, constraints (\ref{battery1}) measure the battery level $b^t$ for each $t \in T$, in agreement with the literature \cite{Cordieri}. At the beginning of the time horizon considered the energy stored in the battery and the energy withdrawn from it are both set to zero by constraints 
(\ref{battery start}) and (\ref{battery out start}). The presence of variables $y^t_{in}$ and $y^t_{out}$ in constraints (\ref{battery out/in}) guarantee that there is no simultaneous withdrawal and injection happening in the battery for every time step $t \in T$. In fact, $y^t_{in}$ and $y^t_{out}$ represent binary variables that take the value 1 if energy is respectively injected in or withdrawn from the battery in $t \in T$ and 0 otherwise. Constraints (\ref{battery out}) and (\ref{battery in}) link the binary variables $y^t_{in}$ and $y^t_{out}$ to the respective continuous variables $b^t_{in}$ and $b^t_{out}$.
 \begin{align}
     \label{battery1}
  b^t = b^{t-1}+\eta_{b}b^t_{in}-\frac{b^t_{out}}{\eta_{b}} \ \forall t \in T\\
  \label{battery start}b^0=0\\
  \label{battery out start}b_{out}^0=0\\ 
  \label{battery out/in}y_{in}^t=1-y_{out}^t \ \forall t \in \{0,..,T\}\\ 
  \label{battery out}b_{out}^t \leq b_{max} y_{out}^t \ \forall t \in \{0,..,T\}\\  
  \label{battery in}b_{in}^t \leq b_{max} y_{in}^t \ \forall t \in \{0,..,T\}
\end{align}
The following set of constraints accounts for the degradation of the battery in time. Constraints (\ref{degradation1}) computes the degradation factor in each $t \in T$, in agreement with the literature \cite{BORDIN2019100824}, which is then applied to the maximum capacity limit in constraints (\ref{maxcapacity}). Finally, constraints (\ref{grid in 3})-(\ref{grid out limit}) limit the energy withdrawn/injected every $t \in T$.   
\begin{align}
  \label{degradation1}
  \frac{B^{fade}}{B^{throughput}}|b^t-b^{t-1}|\leq d^t \ \forall t \in T\\
  \label{maxcapacity}
  b_{max}^t\leq d^t b_{max} \ \forall t \in T 
  \\ 
  \label{grid in 3}
  b_{min}\leq b^t_{in}\leq b_{max}^t \ \forall t \in T\\
  \label{grid out limit}
  b_{min}\leq b^t_{out}\leq b_{max}^t \ \forall t \in T  
\end{align}
Finally, constraints (\ref{definition1})-(\ref{definition3}) define the variables.
\begin{align}
    \label{definition1}
x^t,z^t,b^t,b^t_{in},b^t_{out},e^t_{in},e^t_{out},q^t_{in},q_{solar}^t,m_{ORC}^t \ \geq 0 \ \forall t \in T\\
    \label{definition2}
    g^t \in \mathbb{R} \ \forall t \in T\\
    \label{definition3}
    y_{in}^t,y_{out}^t \in \{0,1\} \ \forall t \in T
\end{align}
\subsection{TET model}\label{TET model}
In the following, an optimization model for Transactive Energy Trading is presented, called the TET model, in a P2P context for a Solar-ORC coupled with a storage system. The TET model is implemented after the S-ORC model to optimize the P2P trading among microgrids. In this sense, a set $N$ of participants in the trading is defined. The participants are essentially consumers or prosumers that can participate as part of the demand, as part of the providers or both for every time step $t$ in which the time range $T$ is divided. The S-ORC model is solved in parallel for every participant, calculating the imbalances produced by over/underproduction of each system. These imbalances were handled in the S-ORC model by the grid, through variables $e_{out}^t$ and $e_{in}^t$. The optimal values of these variables produced by the S-ORC model are then optimized by the TET model. In fact, they are used as parameters $e_{j,out}^t$ and $e_{j,in}^t$. More precisely $e_{j,in}^t$ is the energy needed by participant $j \in N$, while $e_{j,out}^t$ that can be traded by participant $j \in N$. The variables that represent the fluxes of energy that move among participants $j \in N$ every time step $t \in T$ are computed by variables $f_{ij}^t$. The grid is used once again to deal with imbalances. However, this time such imbalances concern the whole system and not the single participant. The pseudo-code of this process are shown in Algorithm \ref{alg:pseudo-code}.\\

\begin{algorithm2e}  
\linespread{0.6}\selectfont
\small
	\KwIn{Set of participants $N$ and set of time steps $T$;}
    \For{$j \in N$} 
    {\KwIn{Parameters for the S-ORC model;}
    Solve the S-ORC model;\\ \KwOut{$e_{j,out}^t$ and $e_{j,in}^t$}}
    \KwIn{For every $j \in N \ \text{and} \ t \in T$ $e_{j,out}^t$ and $e_{j,in}^t$;}
    Solve the TET model;\\
    \KwOut{Optimal solution for the TET model} 
    
	\caption{Pseudo-code of the solving procedure}
	\label{alg:pseudo-code}
\end{algorithm2e}
\vspace{1cm}
The methodology described is functional to the problem we are trying to solve. In fact, we want each single prosumer to first fulfill their own demand, and then to think about trading of residual capacity. Therefore, it seemed more practical to avoid a single optimization model that contained both the TET model and the S-ORC model.\\

The objective of the TET model (\ref{objective function1}) aims to minimize the overall costs of the trading system. Such costs are represented mainly by transmission costs given by parameter $c_T^t$ for every kWh that goes from participant $i \in N$ to participant $j \in N$ every $t \in T$. 
\begin{align}\label{objective function1}
 [\text{TET model}]\    \text{minimize}   \sum_{t \in T}\sum_{i \in N} \sum_{j \in N} c_{T,ij}^t|f_{ij}^t|
\end{align}
Constraints (\ref{flux1}) and (\ref{flux2}) balance the energy sold to another participant or bought from another participant for every participant $i \in N$ and $j \in N$, every time step $t \in T$.
\begin{align}
  \label{flux1}
  \sum_{j \in N} f_{ij}^t\geq e_{j,out}^t \ \forall t \in T,,\ \forall i \in N\\
  \label{flux2}
  \sum_{i \in N} f_{ij}^t\geq e_{j,in}^t \ \forall t \in T,,\ \forall j \in N
\end{align}
Constraints (\ref{flux4}) balance the overall system with the grid, to control over/under production.
\begin{align} 
  \label{flux4}
  h_{in}^t+ \sum_{i \in N} \sum_{j \in N}f_{ij}^t = h_{out}^t \ \forall t \in T  
\end{align}
Constraints (\ref{flux3}) limit the fluxes between all the participants $i,j \in N$ to a minimum and a maximum value, i.e $f^{min}_{ij}$ and $f^{max}_{ij}$, every time step $t \in T$.
\begin{align}
  \label{flux3}
  f^{min}_{ij}\leq f_{ij}^t\leq f^{max}_{ij} \ \forall t \in T, \ \forall i \in N, \forall j \in N
\end{align}
Finally constraints (\ref{final1}) and (\ref{final2}) define the variables.
\begin{align}
     \label{final1}
    f_{ij}^t \in \mathbb{R}\ \forall i,j \in N \ \forall t \in T\\
         \label{final2}
    h_{in}^t, h_{out}^t \geq 0 \ \forall t \in T
\end{align}

\section{Computational experiments}\label{sec:Computational experiments}
In this section, the results obtained by the computational experiments are shown. The computational experiments have been done with four threads with 8 GB, on a computer having
4 cores and a processor Intel(R) Core(TM) i5-7200U @2.50 GHz. All the tests were performed 
using Gurobi 9.1.2  \cite{Gurobi} as a solver. The models were implemented using the JuMP package of Julia \cite{JULIA}.\\

\subsection{Sensitivity analysis on the S-ORC model}\label{sec:Sensitivity analysis}
In this Section, we present the results obtained by performing a sensitivity analysis on the S-ORC model. This was done to inspect, how changing some specifics of the Solar-ORC would affect the system. The demand is represented by an industrial plant, that uses the Solar-ORC in a self-consumption setting. We solve every instance considering a time horizon of one week, with hourly intervals.
All the instances presented were solved within 0.06 seconds.
First, we tested the model using 9 different types of working fluids for the Solar-ORC, to detect the effects that this might have on the mass flow rate. The working fluids have different specifics that are shown in Table \ref{tab:Fluids}.

\begin{figure}  
    \centering
    \includegraphics{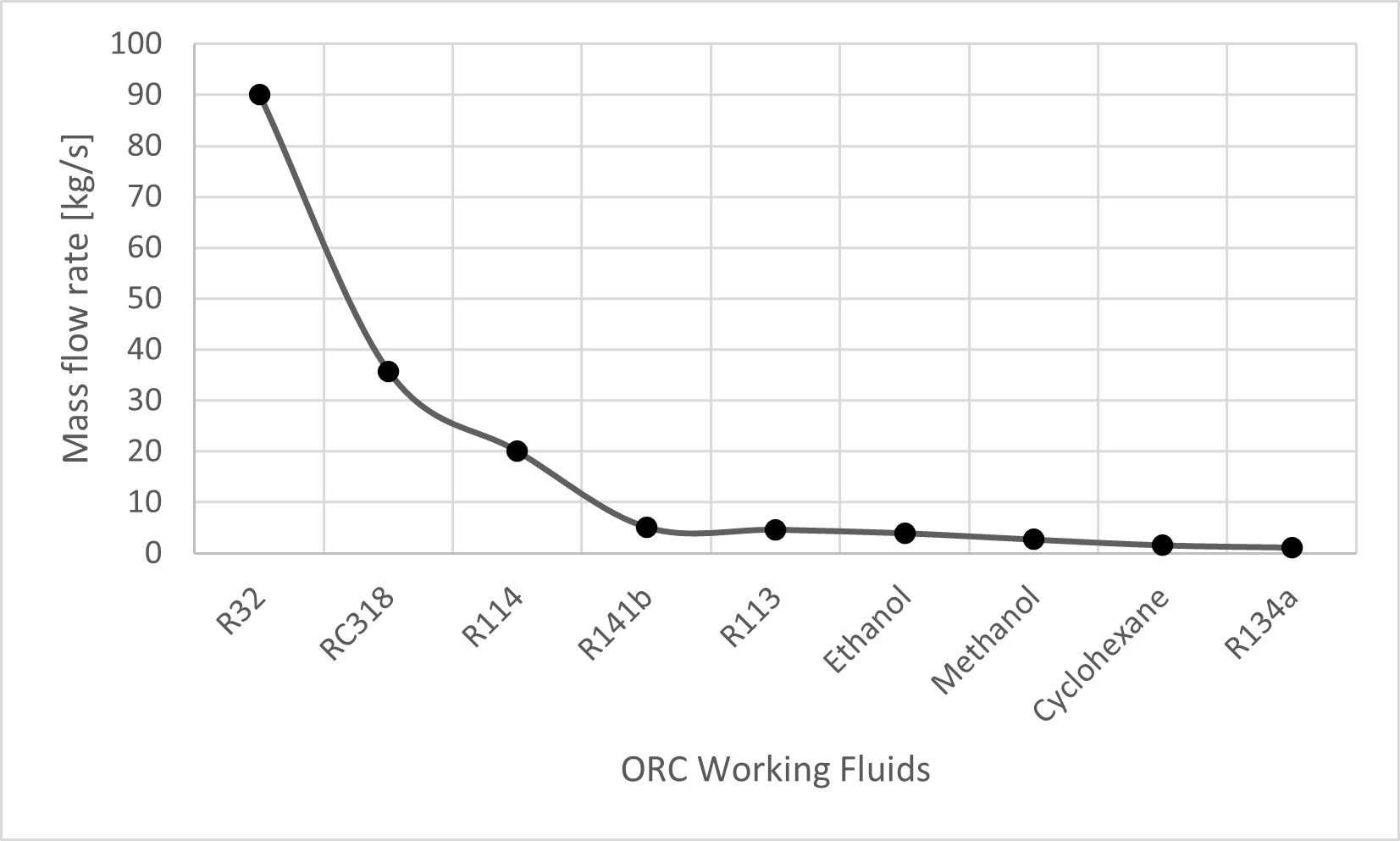}
    \caption{Mass flow rate for different types of working fluids}
    \label{fig:fluids}
\end{figure}

\begin{table}  
    \centering
    \begin{tabular}{|l|l|l|l|l|l|}
    \hline
        Fluid & Molecular Weight [kg/mol] & $T_{crit}$[°C] & $P_{crit}$[MPa] & $C_p$ [J/kg °C] & Density [kg/m\^3] \\ \hline
        Ethanol & 0.046 & 240.8 & 6.148 & 2432 & 0.253100481 \\ \hline
        Methanol & 0.032 & 240.2 & 8.104 & 2512 & 0.369822485 \\ \hline
        Cyclohexane & 0.084 & 280.5 & 4.075 & 154.37 & 0.632911392 \\ \hline
        R134a & 0.102 & 101 & 4.059 & 1268 & 0.8838 \\ \hline
        R141b & 0.11695 & 204.2 & 4.249 & 895 & 0.195 \\ \hline
        RC318 & 0.2 & 115.2 & 2.778 & 898 & 0.028 \\ \hline
        R114 & 0.17 & 145.7 & 3.289 & 845 & 0.05 \\ \hline
        R113 & 0.187 & 214.1 & 3.439 & 867 & 0.215 \\ \hline
        R32 & 0.052 & 78.11 & 5.784 & 848 & 0.011 \\ \hline
    \end{tabular}
    \caption{Working fluids thermodynamic properties}
    \label{tab:Fluids}
\end{table}
Figure \ref{fig:fluids} shows the mass flow rate needed using different fluids for a 2 kW Solar-ORC, considering the same working conditions. As one can observe working fluids like Ethanol, Methanol, Cyclohexane, and R134a need a lower mass flow rate. Thus, from an economic perspective, these fluids can be interesting, especially for large-capacity systems.

Subsequently, we analyzed the consequences on the system's performance by considering 9 different sizes of the ORC, shown in Table \ref{tab:ORC}.

\begin{table}  
    \centering
    \begin{tabular}{|l|l|}
    \hline
        Size [kW] \\ \hline
        0.5\\ \hline
        1 \\ \hline
        1.5 \\ \hline
        2 \\ \hline
        2.5 \\ \hline
        3 \\ \hline
        3.5 \\ \hline
        4 \\ \hline
        4.5 \\ \hline
    \end{tabular}
        \caption{Organic Rankine Cycle sizes}
    \label{tab:ORC}
\end{table}

\begin{figure}  
    \centering
    \includegraphics{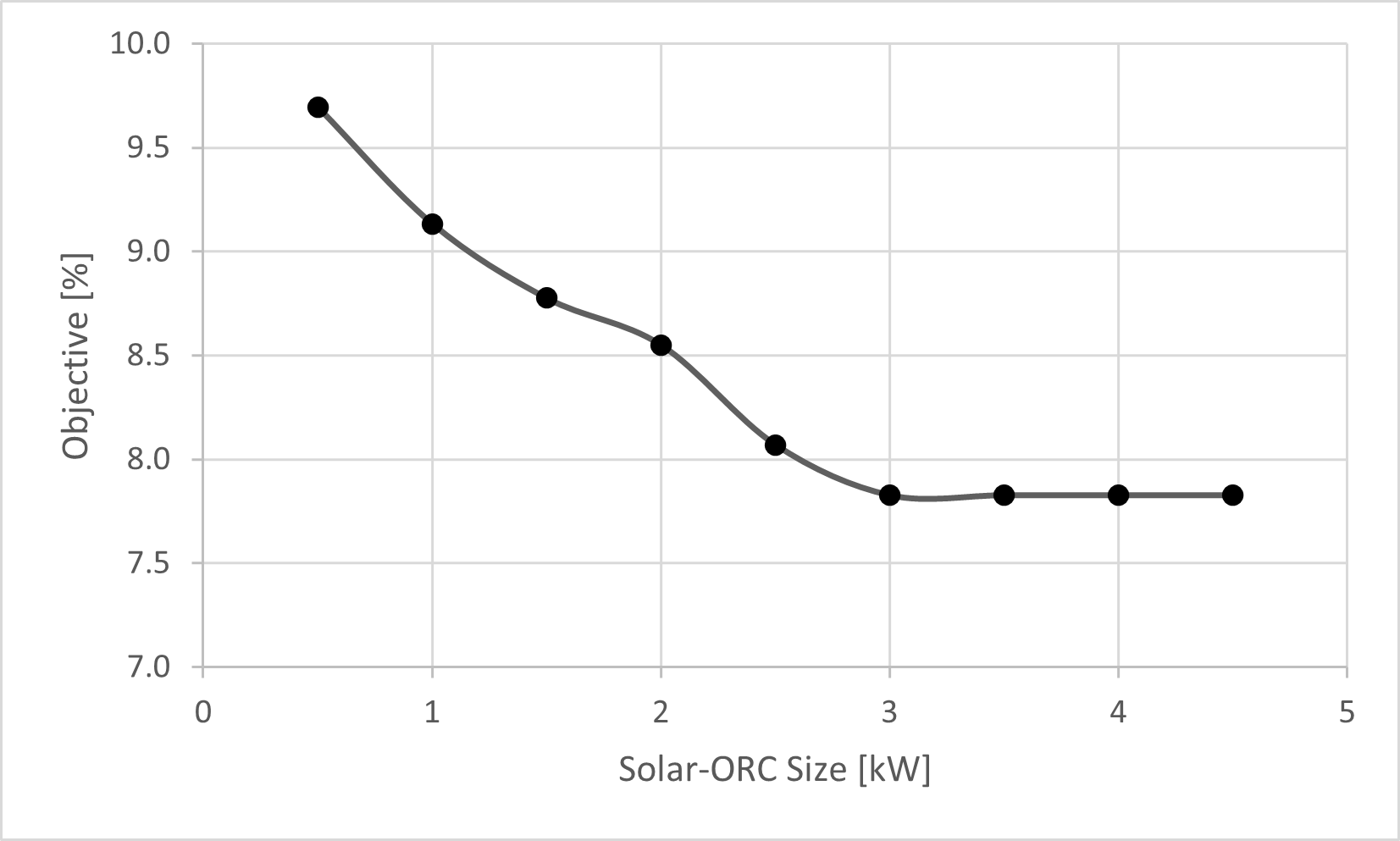}
    \caption{Objective difference with the size of Solar-ORC}
    \label{fig:sizes results}
\end{figure}

Figure \ref{fig:sizes results} shows the difference in objective function considering Ethanol as a working fluid. As one can observe, the objective decreases by increasing the size of the plant, up to a certain threshold. The decrease in the objective is due to a decrease in the electricity provided by the grid, to satisfy the final consumer's demand. When the optimal size to satisfy such demand is reached, there is no economic benefit to increase the plant's size further. This is consistent with the self-consumption framework that we are considering. 

Then, we studied the impact of the solar collector's efficiency on the system. 5 different ORC and solar collector couplings were tested. The specifics regarding the solar collectors are shown in Table \ref{tab:solar collectors}.  The tests were performed considering the same conditions for all the solar collectors, using literature values \cite{KALOGIROU2004231}.

\begin{table}  
    \centering
    \begin{tabular}{|l|l|}
    \hline
        Solar collector & Efficiency [\%]\\ \hline
        FPC & 65\\ \hline
        ETC & 87\\ \hline
        CPC & 65\\ \hline
        PTC & 85\\ \hline
        LFR & 66\\ \hline
    \end{tabular}
        \caption{Solar collector's specifics}
    \label{tab:solar collectors}
\end{table}

\begin{figure}  
    \centering
    \includegraphics{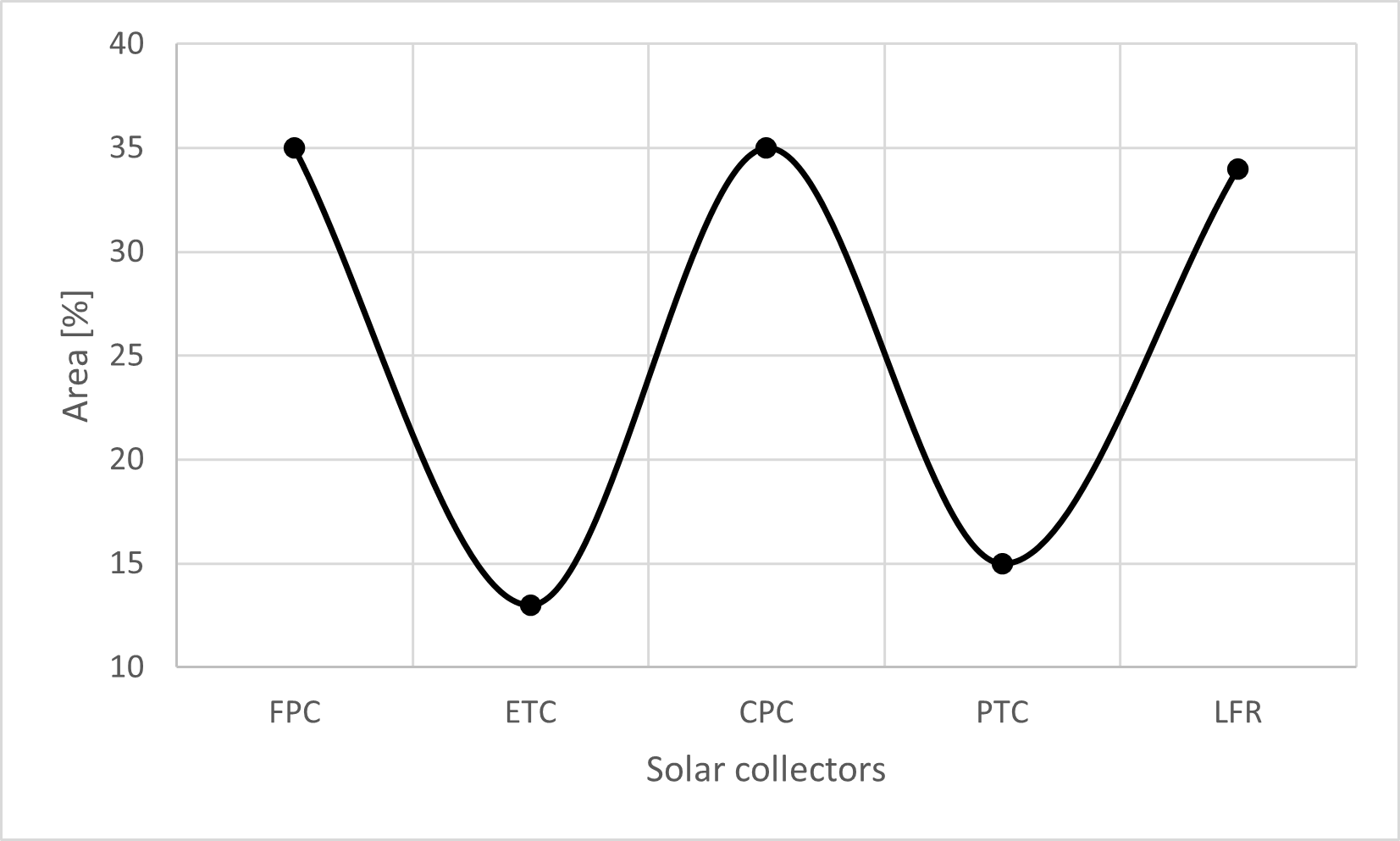}
    \caption{Solar-ORC size difference changing the solar collector's technology}
    \label{fig:solar collectors}
\end{figure}

The main value impacted by the change in the solar collector's efficiency was the sizing of the ORC in mass flow rate. This impact is shown in Figure \ref{fig:solar collectors}.

As one can observe in Figure \ref{fig:solar collectors} ETC and PTC technologies seem more favorable for this kind of application.\\

\begin{figure}  
    \centering
    \includegraphics{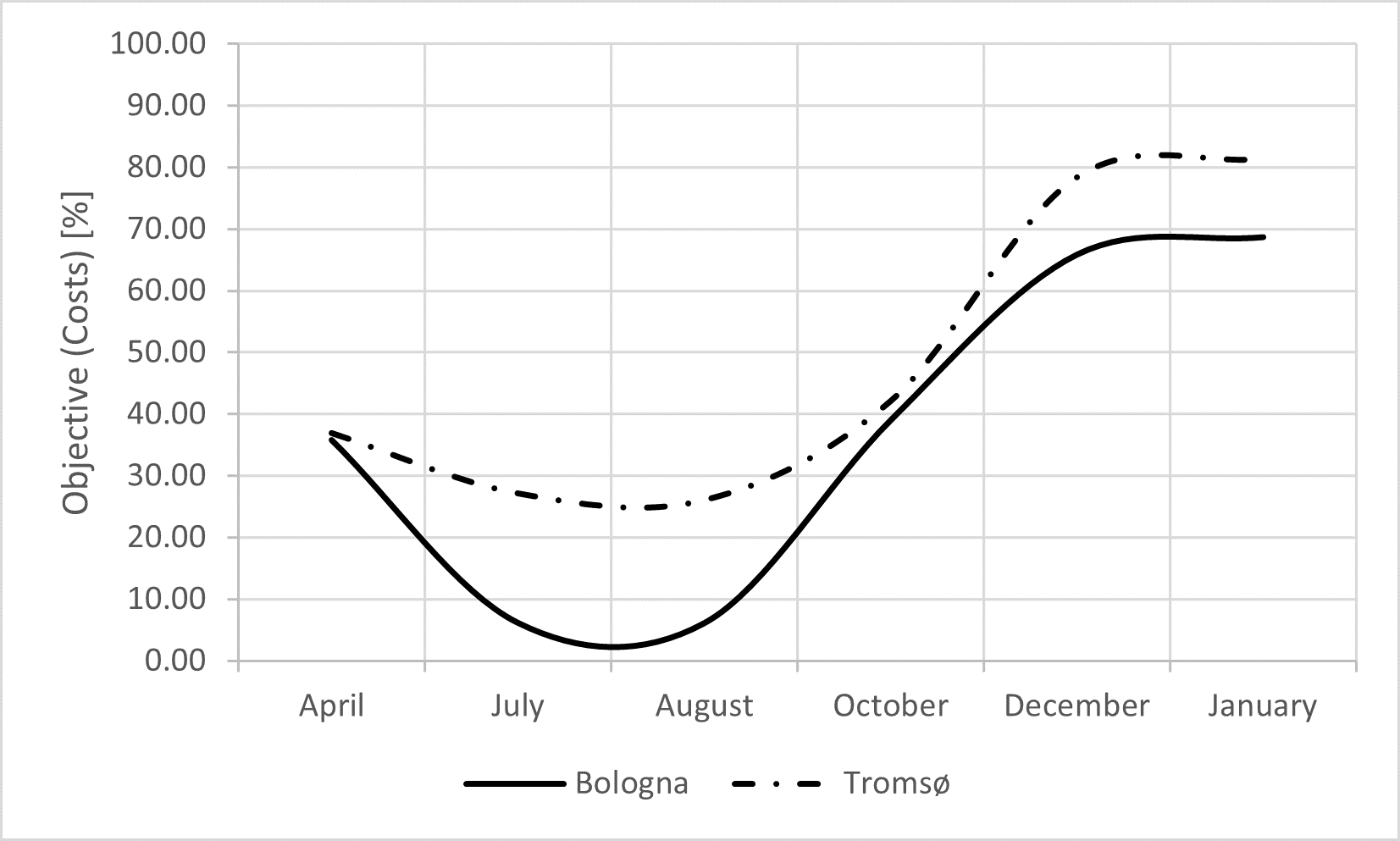}
    \caption{Comparison of objectives for two different locations for six significant weeks}
    \label{fig:Tromsø vs Bologna1}
\end{figure}

Later we inspected the effect of weather conditions on the system. We detected four representative weeks in the months of April, July, October, and January. Moreover, we considered two locations for the system: the city of Bologna in Italy, and the city of Tromsø in Norway. These cities represent two completely opposite scenarios, that could both potentially benefit from the system considered. In fact, Bologna (44.4949° N, 11.3426° E) is located in the northern part of Italy in the Emilia-Romagna region. It has a typically humid temperate climate with cold, humid winters and hot, muggy summers. Precipitation is moderate, while the rains are fairly well distributed throughout the year, even if two maxima are noted in spring and autumn, and two relative minima in winter and summer. On the other hand, Tromsø (69.6492° N, 18.9553° E) is a city in Northern Norway located in the county of Troms and Finnmark. It is subject to a subarctic climate, with very cold winters and cool summers. Since we are north of the Polar Circle, the sun does not rise (polar night) from November 28th to January 14th, while it does not set (midnight sun) from May 19th to July 26th.\\ 
We present results for a 2 kW Solar-ORC using Ethanol as a working fluid. We make the hypothesis by considering the same cost of electricity sold by the grid both in Bologna and Tromsø. This is done to highlight the real difference in terms of solar incidence between these two locations. This hypothesis stands for all the following tests unless specified.
\begin{figure}  
    \centering
    \includegraphics{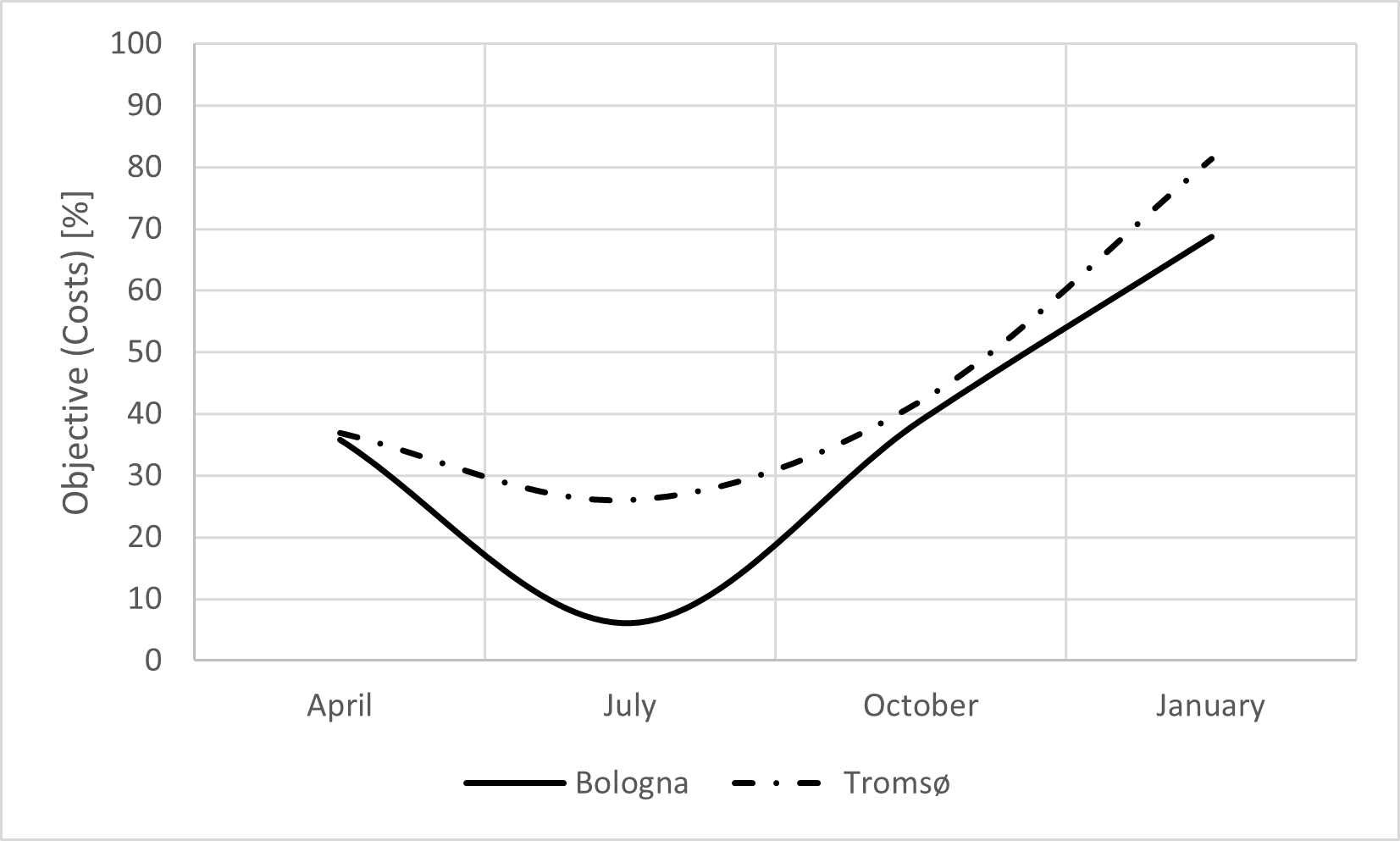}
    \caption{Comparison of objectives for two different locations for four significant weeks}
    \label{fig:Tromsø vs Bologna}
\end{figure}
As shown in figure \ref{fig:Tromsø vs Bologna} Bologna as a location gains a higher economic advantage than Tromsø, regardless of the season. However, this advantage is more significant during July and January. This is consistent with the weather conditions in Tromsø during winter, with almost no solar incidence. Moreover, regardless of the midnight sun phenomenon occurring in Tromsø during summer, the solar incidence in Bologna is still higher due to its latitude.\\
We introduce also two more representative weeks in December and August, as shown in Figure \ref{fig:Tromsø vs Bologna1}. In fact, we want to inspect further possible minimum and maximum points. As one can observe, while the real minimum point is going to be between the last week of July and the first of August for both cities, the real maximum point is more evident in December for Tromsø than it is for Bologna where December and January are almost equivalent.
Finally, we performed a comparison in terms of operational costs between using the Solar-ORC as an electricity source, and using the grid as the only source, obtaining an improvement of 12\% on average, 10\% for Tromsø and 14\% for Bologna.

\subsection{Computational experiments on the S-ORC model}\label{sec:Computational SM model}
\begin{figure}  
    \centering
    \includegraphics[width=\textwidth,keepaspectratio]{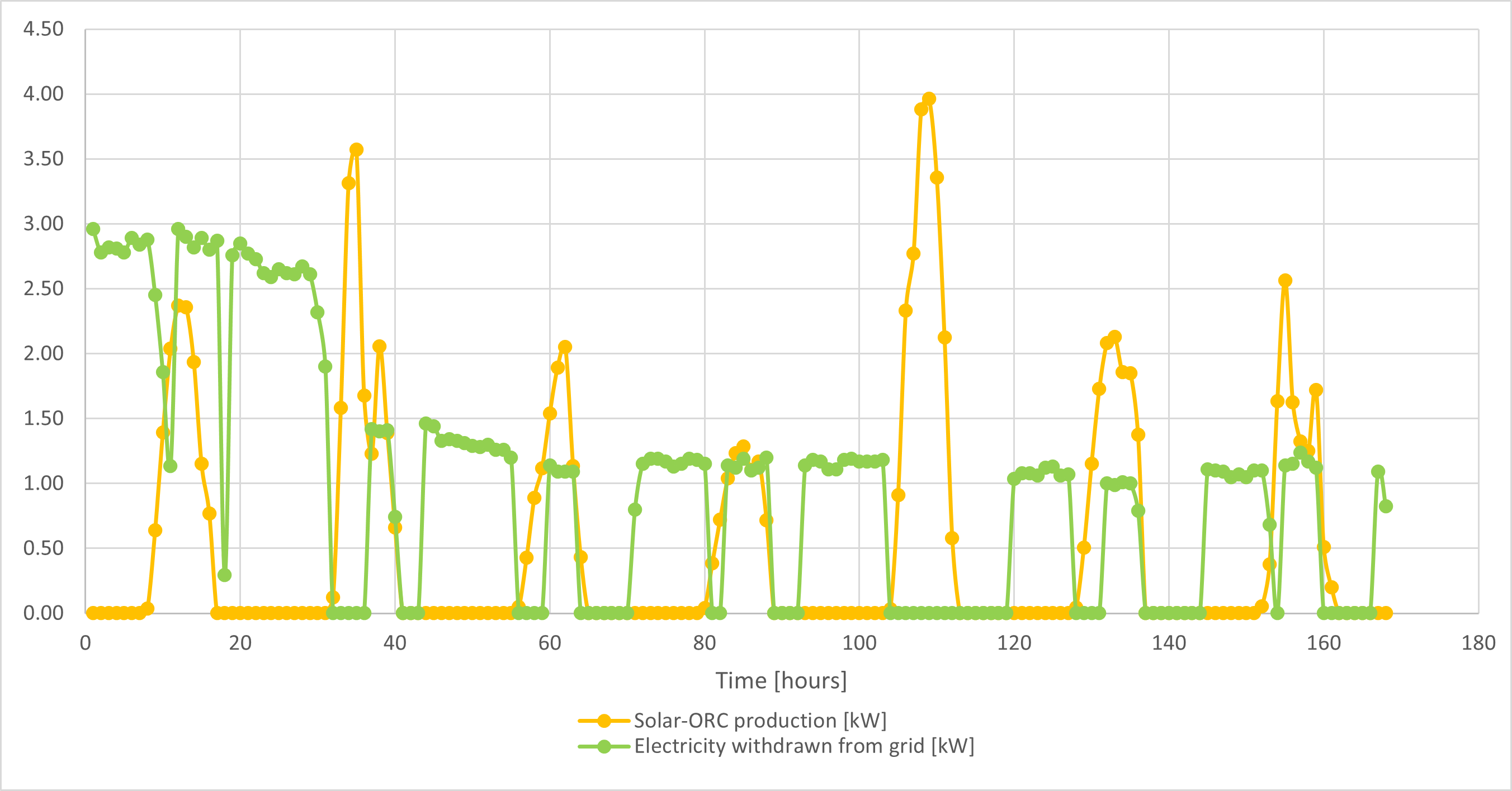}
    \caption{Solar-ORC production and electricity withdrawn from the grid for the city of Bologna}
    \label{fig:production vs grid}
\end{figure}

\begin{figure}  
    \centering
    \includegraphics[width=\textwidth,keepaspectratio]{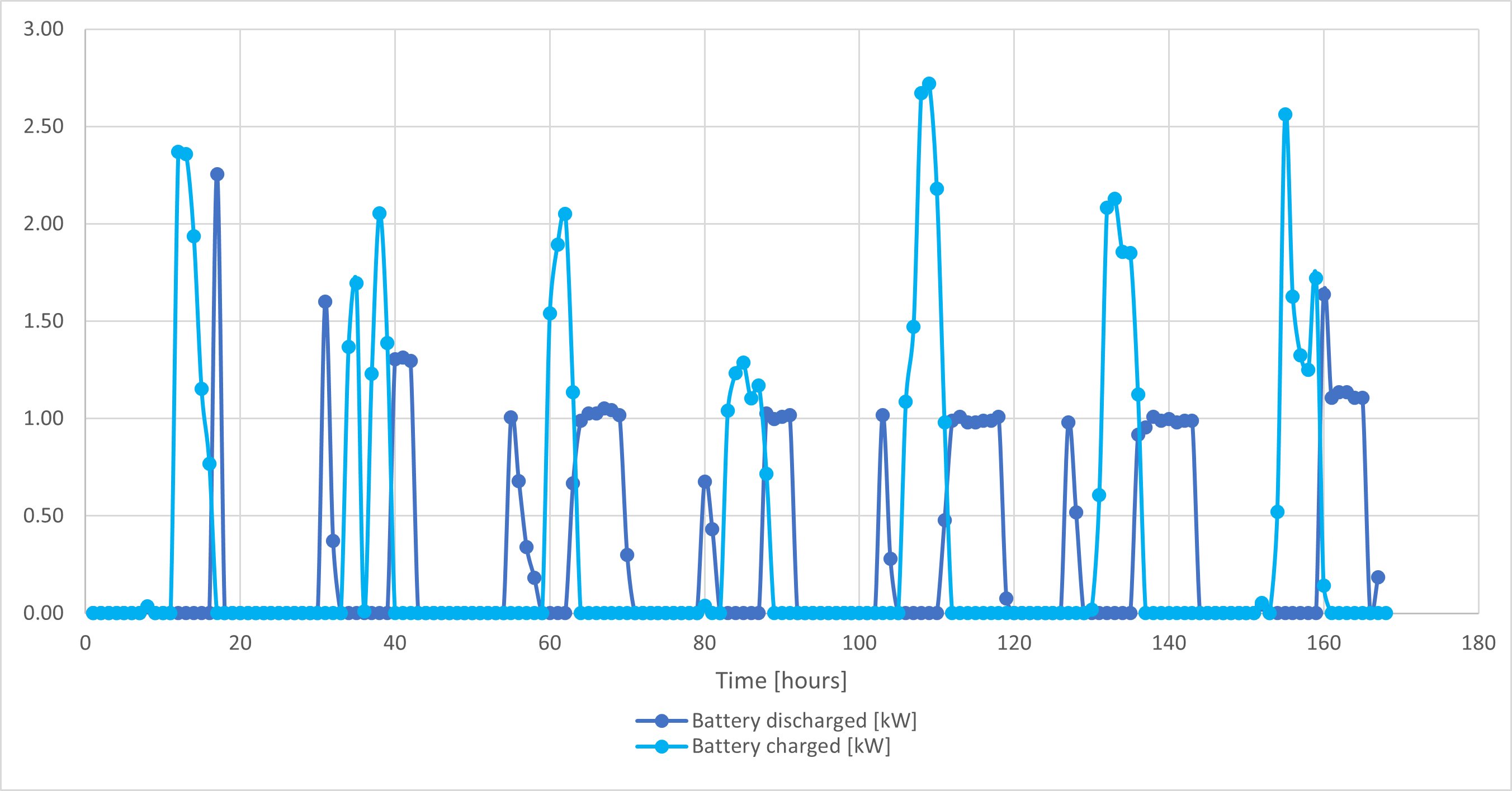}
    \caption{Battery charge and discharge for the city of Bologna}
    \label{fig:battery charge adn discharge}
\end{figure}

This Section presents the results given by testing the S-ORC model. The results shown are referred to the city of Bologna.
Figure \ref{fig:production vs grid} and figure \ref{fig:battery charge adn discharge} show the results given by the model in terms of Solar-ORC production, electricity withdrawn from the grid, and battery usage.

As one can observe in Figure \ref{fig:production vs grid} the grid is used by the S-ORC model to compensate for the lack of production by the Solar-ORC. Specifically, it is used most during the first time periods when there is no solar energy available and the battery is still not charged. This behavior is consistent with the realistic management of the plant. Moreover, as shown in Figure \ref{fig:battery charge adn discharge} usage of the battery is accordingly to its charging operations. In fact, the withdrawal from the battery starts consequently to charging operations, meaning there is actually available energy in the storage system.\\

\subsection{Computational experiments on the TET model}\label{sec:Computational TET model}
\begin{figure}  
    \centering
    \includegraphics[width=\textwidth,keepaspectratio]{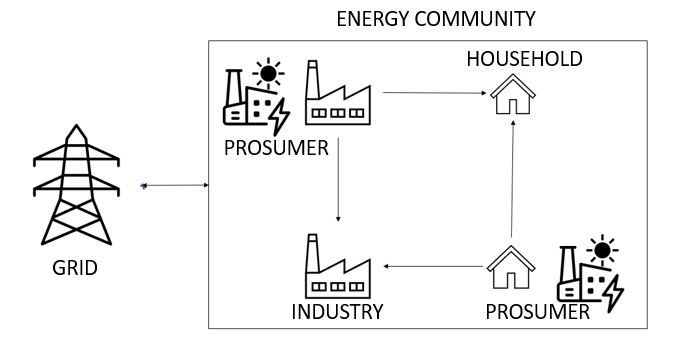}
    \caption{Sytem with different types of consumers, some of them being prosumers supplied by a Solar-ORC}
    \label{fig:system}
\end{figure}

\begin{figure}  
    \centering
    \includegraphics[width=\textwidth,keepaspectratio]{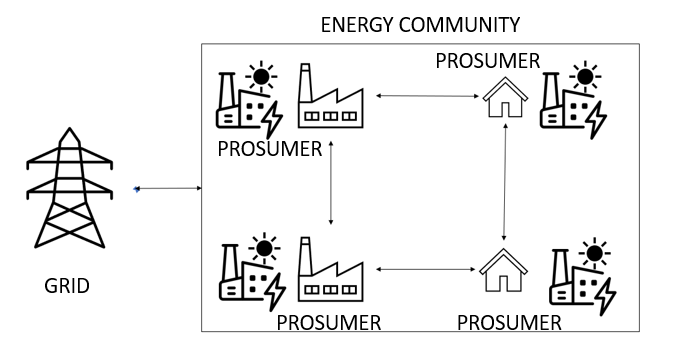}
    \caption{Sytem where all consumers are prosumers supplied by a the Solar-ORC}
    \label{fig:system 1}
\end{figure}

\begin{figure}  
    \centering
    \includegraphics{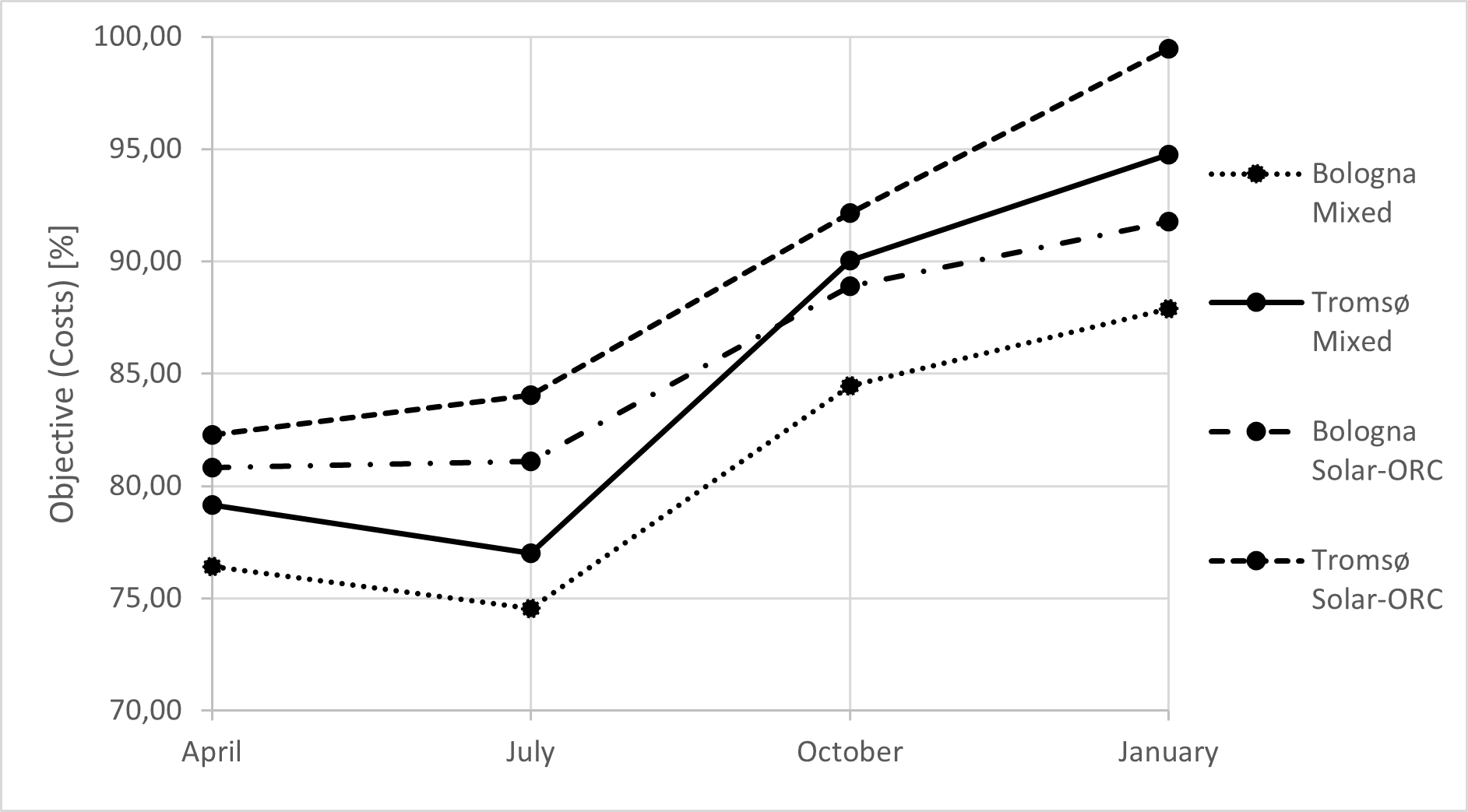}
    \caption{Results from the TET model tested on different instances}
    \label{fig:results total}
\end{figure}

A second step was to inspect the role that the Solar-ORC can potentially have when introduced in a peer-to-peer context. The TET model has been tested with multiple instances. Each instance represents a different energy community. At first, we consider an instance that is described in Figure \ref{fig:system}. Here, the components of the system are partly consumers and partly prosumers. The prosumers are supplied energy by a Solar-ORC in a self-consumption framework, thus the Solar-ORC should first satisfy their demand and then the other consumers' demand. We introduce different types of demand, i.e. industries (that work 24/7) and households. 
Then, we consider an instance where all the consumers become prosumers, each of them being supplied by a different Solar-ORC. The system of the instance is represented in Figure \ref{fig:system 1}.  
The results of such instances are shown in Figure \ref{fig:results total}. 
The introduction of Solar-ORC in a P2P trading context results in an economic gain of an average of 4\%, compared to having a single Solar-ORC. This average is measured as a mean of the results coming from considering both the schemes shown in Figure \ref{fig:system} and Figure \ref{fig:system 1}. The tests show a higher gain in terms of total operational costs for the Bologna location, an average of 4.5\%. This is coherent with the results obtained by the S-ORC model, as shown previously. The results obtained in Bologna by the S-ORC model are still better in terms of operational costs than the ones obtained in Tromsø by the TET model. This is justified by the latitude of Bologna and its greater solar irradiation.  
 Moreover, the gain increased during the summer and spring seasons, in both locations. This is also consistent with the higher solar activity of such periods, as previously discussed.\\
 Finally, we performed an overall comparison in terms of operational costs between performing TET using the Solar-ORC as an electricity source, and using the grid as the only source, and the improvement was around 16\% on average, 19\% for Bologna and 14.7\% for Tromsø.

\section{Interpretation of the results}\label{Discussions}
This section discusses further the results presented in section \ref{sec:Computational experiments}. The results show that both the S-ORC model and the TET model produce realistic solutions. The use of Solar-ORC resulted in an improvement of 12\% on average ( 10\% for Tromsø and 14\% for Bologna) in terms of operational costs, compared to not using such technology. Moreover, when introduced in a P2P trading context the improvement is even greater, around 16\% on average (19\% for Bologna and 14.7\% for Tromsø).
These solutions change reasonably with the weather conditions both in terms of latitude and season. In fact, given the same cost of electricity, Bologna tends to have greater economic gain compared to Tromsø. The difference varies between 1\% and 20\%, with a medium value of 11.72\%. This is consistent with the difference in solar incidence that these two locations have, because of their latitude. Furthermore, this difference is enhanced during the winter season, especially in December and January, where the difference in terms of costs is around 20\%. As a matter of fact, during these months Tromsø is subject to the "Polar Night" phenomenon, with almost 24 hours of darkness. On the contrary, during springtime, the difference in the objective is almost none, around 1\%, consistently with the weather in Bologna being often cloudy.\\
Concentrating just on the Tromsø instance, one can observe that the results change when taking into account local electricity prices. In fact, the objective values are drastically decreased due to lower local electricity prices, especially during summer when they are extremely low. Given that the electricity prices in Norway are significantly lower than in Italy, especially during summer, it might seem not convenient to invest in a renewable energy source such as Solar-ORC. However, recently new studies, i.e. Nguyen et. al \cite{9967828}, have concentrated on solar power systems in this area, showing how Tromsø may profit from them. According to Eikeland et. al \cite{9967829}, the energy production of these systems could be coupled with Cruise ships' energy consumption.
Cruise ships have a great demand, that needs to be fulfilled even when they are located in a harbor. If the energy supply from the harbor is not sufficient, the ships need to run their own motors to operate, thus producing a great environmental impact. Furthermore, Eikeland et. al \cite{9967829} state that the highest number of visiting Cruise ships (C.S) is during the tourist season in June, July, and August. This period coincides with the "Midnight Night" phenomenon in Tromsø, with almost 24 hours of solar power availability.\\
Another way to increase the exploitability of this technology in Artic areas like Tromsø could be the coupling with seasonal storage. Considering the aforementioned peculiar climate of this area, long-term storage seems even more appropriate. This way, the great solar energy availability in summer would be exploited also during winter. Thus, the economic advantage would be evenly distributed throughout the year.\\

Both the S-ORC model and the TET model solve all the instances within a few seconds. This is reasonable with the weekly time range of optimization. This paper focuses on this time range, as it is more relevant to the short-to-medium-term operational planning we are examining. Moreover, the TET model solves first multiple S-ORC models in parallel (one for each micro-grid) and then optimizes the trading within the S-ORC models. Thus, the computational time is given by the addition of the longest computational time among the parallel S-ORC models and the trading part computational time. This solution is less time costly than considering a single general optimization of the system. In fact, the variables that would be handled simultaneously in one single model by the solver, are split into different models and thus into different computational procedures. Moreover, having different computational procedures is more consistent with the self-consumption framework we are considering. In fact, in a self-consumption framework, the final user's goal is first to fulfill its demand through its own power production plant and then to adjust over/undersupply.  
In view of this concept, the goal of our model in the trading phase is to only optimize the under/oversupply of electricity among the micro-grids of our system.\\
Being computationally tractable, these modeling approaches can be of great value for a wide variety of energy and power systems tools. Indeed, they can be introduced in larger energy and power systems models, if the user wishes to investigate solar ORC in a broader context. Open source tools such as PyPSA, \cite{PyPSA}, highRES \cite{highRES}, GenX from MIT \cite{GenX}, Sienna from NRL \cite{SIENNA}, could benefit from such a technology-oriented approach to be included, for instance, as a module for more specialized studies.\\

\subsection{Value of the results for solar power plants research}\label{sec:value solar}
In this work, we presented a peculiar application of solar energy, to reduce the emission impact of conventional power production cycles.  
The tests presented in section \ref{sec:Computational experiments} show a promising applicability of solar-ORC in peer-to-peer transactive energy trading. Therefore, this technology can be a profitable application of solar energy, opening a new perspective in this field.

The algorithm presented results to be an effective tool to support the operational management of such a system, reaching optimal solutions. Moreover, it is able to perform for instances of weekly time horizons, with competitive computational time. This opens the possibility for the use of such an algorithm in realistic instances, as a support tool for decision-making actors. 

The proposed formulation can also be integrated within wider energy and power systems optimization models to explore new synergies and scenarios. This applies, for example, to the integration between solar ORC and pumped thermal electricity storage \cite{frate2022assessment}. Models such as \cite{ibrahim2017pumped} can be further expanded by including a technological solar ORC module to investigate the impact on active distribution network applications. Other topical optimization problems such as the reliability-oriented network restructuring problem defined in \cite{mishra2019rnr} and \cite{bordin2021multihorizon} could greatly benefit from a solar ORC module to enhance analyses and investigate new local power mitigation strategies in renewable settings.

ORC mathematical optimization can also play a role in green transportation \cite{bordin2021behavioural} by efficiently converting low-grade waste heat into usable electricity.
This technology enhances energy efficiency in electric vehicle (EV) charging systems, reduces costs, and supports the use of renewable energy \cite{imran2019optimization}. By optimizing energy recovery, ORC systems help minimize environmental impact, contributing to the decarbonization of the transportation sector and promoting more sustainable mobility solutions.

ORC optimization enhances wind energy systems by improving efficiency and stabilizing grid integration. It converts waste heat into additional electricity, balancing supply and demand while storing excess energy \cite{nematollahi2019evaluation}. When integrated with predictive models \cite{mishra2020features} for wind ramp events, the ORC model improves decision-making by anticipating fluctuations, thus supporting grid stability and maximizing the value of renewable resources.

The examples and reflections above demonstrate how the proposed ORC optimization model can be seamlessly integrated into a wide range of existing predictive and prescriptive analytics tools. By enhancing both energy system modeling and decision-making processes, it enables more comprehensive analyses, fosters synergies across various renewable energy applications, and supports the development of innovative strategies for improving efficiency, sustainability, and grid stability.

In the sensitivity analysis proposed in section \ref{sec:Sensitivity analysis} we inspected the role of the solar collector's technology. The outcome of such analysis confirmed the impact of the solar collector, but mainly on the sizing of the Solar-ORC. In fact, the results highlight an impact on the dimensions of the Solar-ORC. This implies an influence mainly on the design phase, thus on the capital costs' definition. However, the aim of our work is to give an insight into operational management, thus on factors that may affect the operational functioning of the system. Therefore, we do not consider capital costs and the design process. The fundamental idea is that the system proposed has already been designed and installed. 

We can conclude then, that the typology of this component can be neglected, as far as operational decisions are concerned. In fact, once optimally projected, its capacity is optimized to be compliant with the system's needs.

\subsection{Limitations of the work and future research}\label{Future work}
The results shown in section \ref{sec:Computational experiments} highlight an overall gain in the implementation of a predictable and manageable system as the one we present in this paper for a P2P transactive energy trading context. An investment costs analysis is not included in the model that is focusing on the operational aspects. Having investment costs would give for sure a more complete view of the real gain of such systems. However, the scope of this paper was to produce an optimization model for the operational management of the system. The main assumption is that all the technologies involved are already installed. In fact, our purpose was not to discuss strategies connected to the investment phase. Therefore, a suggestion for future directions related to this work is to concentrate on the integration of investment decisions. An optimization model for the design part of the solar-ORC can be included and solved before the operational management optimization problem is discussed in this paper. The methodology would then include three steps: a design optimization for every single Solar-ORC, an operational management optimization of every single Solar-ORC, and finally a TET optimization of the local community.

The role of seasonal storage could also be a future direction to investigate. Seasonal storage is a highly discussed topic, especially in the hydrogen field. Technologies like Power to hydrogen (P2H) or Power to gas (P2G) are able to produce hydrogen during high peak production periods of RES and store it for long periods. For peculiar locations like Tromsø, where the weather conditions change drastically with the season, the addition of seasonal storage could improve the final result. This way, the benefits deriving from RES production could be extended to the whole year. It would be interesting to see how the presence of seasonal storage could impact the efficiency of the system analyzed in this paper and the economic gain for the Tromsø instance. In this case, it would be necessary to modify the SM model to include constraints regarding seasonal storage.  

Finally, further analysis could be made focusing on the impact of the real cost of electricity for the Tromsø instance. As we mentioned in section \ref{Discussions}, in this study we concentrate on the meteorological impact on the system. For this reason, the role of electricity costs is not inspected, thus we use the same cost of electricity for Bologna and Tromsø. More specifically, we apply the cost of electricity given by the Italian electricity market both for Bologna and Tromsø. This could lead to an over/underestimation of the real benefit for the Tromsø instance. Thus, future works could inspect the role of electricity costs. This could be easily done by using the real costs of electricity in Tromsø.

{\section{Conclusions}\label{sec:Conclusions}
This paper investigates the potential that a power generation technology like a Solar-ORC could have in being introduced in a P2P transactive energy trading context.
In sight of this, a tool based on operation research techniques has been implemented, able to optimize the scheduling of both the Solar-ORC and the trading process. First, a MILP model for the operations scheduling of the Solar-ORC has been developed, the S-ORC model. Then, a MILP model has been formulated for the P2P transactive energy trading between multiple prosumers in a local energy market where some Solar-ORCs are present as power generation plants owned by some prosumers, the TET model.

First, a sensitivity analysis on the S-ORC model is performed, to inspect, how changing some specifics of the Solar-ORC would affect the system. We inspected the effects on the system given by different types of working fluids, different sizes of the Solar-ORC, and different weather conditions. Every instance has been tested for hourly intervals within a time horizon of one week.
Each instance presented was solved within 0.06 seconds.
From an economic perspective,  fluids like Ethanol, Methanol, Cyclohexane, and R134a are potentially more valuable, especially for large-capacity systems. In fact, they need a lower flow rate compared to others, given the same weather conditions and size of the plant.

On the contrary, when considering the same working fluid and weather conditions, the objective decreases by increasing the size of the plant, up to a certain threshold. In fact,  when the optimal size to satisfy the demand of the prosumer is reached, there is no economic benefit to increasing the plant's size further. This is consistent with the self-consumption framework that we are considering.

Later, five different ORC and solar collector couplings have been tested, to observe the impact of the solar collector on the system. Specifically, we used the solar collector's efficiency as the main parameter. The main value that was impacted by this parameter was the mass flow rate of the ORC cycle. This reflects a change mainly in the sizing of the cycle, thus an adjustment of design parameters of the system, not on operational management.

We then inspected the effect of weather conditions on the system, considering six representative weeks taken from different seasons and two locations for the system: the city of Bologna in Italy, and the city of Tromsø in Norway, with the same plant size and working fluid.
Regardless of the season, Bologna as a location gains a higher economic advantage than Tromsø, with two significant differences between July and August, and December.

The proposed S-ORC model closely resembles the real-world production process in the power plant. The grid is used by the S-ORC model to compensate for the lack of production by the Solar-ORC, while the usage of the battery is accordingly to its charging operations.

Later the TET model has been tested with multiple instances. Each instance represents a different energy community. The TET model was able to solve all the instances within a few seconds, giving reasonable results for all the prosumers involved. Coherently to the S-ORC model, the tests show a higher gain in terms of costs for the Bologna location, around 4.5\% in terms of operational costs.

In conclusion, the results highlight an overall gain in the implementation of a predictable and manageable system as the one we present in this paper for a P2P transactive energy trading context, on average 16\% in terms of operational costs. Since the aim of this paper was to produce an optimization model for the operational management of the system, an investment costs analysis is not included. We want to specify that even though the implemented methodology has been tested for two specific cases (Bologna and Tromsø), it can be applied to different case studies. In fact, parameters can be modified to adapt to the specifics of the case study considered. Thus, the model will be made available in the GitHub platform under the name ``OPTI-ORC" 
 (https://github.com/sambeets/OPTI-ORC). Future directions related to this work would be to concentrate on the integration of investment decisions. Moreover, the introduction of long-term storage systems in Arctic areas like Tromsø could be another suggestion for future studies.

\bibliographystyle{IEEEtran}
\bibliography{sn-article}

\end{document}